\documentclass{amsart}

\usepackage[a4paper,top=4cm,bottom=4cm,left=3cm,right=3cm]{geometry}

\usepackage{amssymb, amsmath, amsthm, amsfonts, euscript, color}
\usepackage{tabularx}
\usepackage{array}
\usepackage{graphicx} 
\usepackage{comment}
\usepackage[hypertexnames=false,
    backref=page,
    pdftex,
    pdfpagemode=UseNone,
    breaklinks=true,
    extension=pdf,
    colorlinks=true,
    linkcolor=blue,
    citecolor=blue,
    urlcolor=blue,
]{hyperref}

\theoremstyle{plain}
\newtheorem{theorem}{Theorem}[section]
\newtheorem{lemma}[theorem]{Lemma}

\newtheorem{corollary}[theorem]{Corollary}

\newtheorem*{theorem*}{Theorem}

\theoremstyle{definition}
\newtheorem{definition}[theorem]{Definition}

\newtheorem{example}[theorem]{Example}
\newtheorem*{proposition 4.1}{Proposition 4.1}
\newtheorem*{corollary*}{Corollary}
\newtheorem*{example*}{Example}
\newtheorem*{definition*}{Definition}

\theoremstyle{remark}
\newtheorem{remark}[theorem]{Remark}

\newcommand{\Q}{\mathbb{Q}} 
\newcommand{\R}{\mathbb{R}} 
\newcommand{\C}{\mathbb{C}} 
\newcommand{\del}{\partial}
\newcommand{\delbar}{\overline{\partial}}
\DeclareMathSymbol{\Finv} {\mathord}{AMSb}{"60}

\title{p-K\"ahler structures on compact complex manifolds}
\author{Ettore Lo Giudice}
\address[Ettore Lo Giudice]{
Dipartimento di Scienze Matematiche, Fisiche e Informatiche
Unit\`a di Matematica e Informatica\\
Universit\`a degli Studi di Parma\\
Parco Area delle Scienze 53/A, 43124\\
Parma, Italy}
\email{ettore.logiudice@unipr.it}

\thanks{The author is partially supported by GNSAGA of INdAM}

\begin{document}

\maketitle
\tableofcontents
\begin{abstract}
    Let $(M,J)$ be a complex manifold of complex dimension $n$. A $p$-K\"ahler structure on $(M,J)$ is a real, closed $(p,p)$-transverse form. In this paper, we address the conjecture of L. Alessandrini and G. Bassanelli on $(n-2)$-K\"ahler nilmanifolds equipped with nilpotent complex structures and holomorphically parallelizable nilmanifolds. We also derive necessary conditions for the existence of smooth curves of $p$-K\"ahler structures, starting from a fixed $p$-K\"ahler structure, along a differentiable family of compact complex manifolds. In addition, we study the cohomology classes of $p$-K\"ahler (resp. $p$-symplectic, $p$-pluriclosed) structures on compact complex manifolds. We provide several examples of families of compact complex manifolds admitting $p$-K\"ahler or $p$-symplectic structures.  
\end{abstract} 

\section{Introduction}\label{Section 1}
Let $(M,J)$ be a compact complex manifold of complex dimension $n$. It is well known that the existence of K\"ahler metrics on $(M,J)$ imposes strong topological constraints. For example, if $(M,J)$ admits a K\"ahler metric $g$, whose fundamental form is $\omega$, then the $p$-th power of $\omega$, for $1 \leq p \leq n$, defines a non-zero class in the de Rham cohomology of $M$. Hence, the even index Betti numbers are positive. Moreover, the odd index Betti numbers are even, the complex manifold underlying the K\"ahler manifold satisfies the {\em $\del \delbar$-lemma}, i.e., 
\begin{equation*}
    \text{Ker $\del$ $\cap$ Ker $\delbar$ $\cap$ Im $d$ = Im $\del \delbar$},
\end{equation*}
and the manifold $M$ is {\em formal} in the sense of Sullivan \cite{DGMS1975}. Furthermore, the complex de Rham cohomology of compact K\"ahler manifolds decomposes as the direct sum of Dolbeault cohomology groups. In particular, {\em $(p,q)$-Dolbeault cohomology groups} are isomorphic to {\em $(p,q)$-Bott-Chern cohomology groups} and {\em $(p,q)$-Aeppli cohomology groups}. The definitions of Aeppli and Bott-Chern cohomology groups are the following
\begin{equation*}
    H^{p, q}_{BC}(M) = \frac{\text{Ker} \, \del \cap \text{Ker} \, \delbar}{\text{Im} \, \del \delbar}, \quad H^{p, q}_{A}(M) = \frac{\text{Ker} \, \del \delbar}{\text{Im} \, \del \oplus \text{Im} \, \delbar},
\end{equation*}
where the operators act in the correct space of forms. Moreover, by a remarkable result of K. Kodaira and D.C. Spencer \cite{KS1960}, small deformations of compact K\"ahler manifolds are still K\"ahler manifolds. 

However, several complex manifolds fail to admit K\"ahler metrics. Despite this, they may still support special Hermitian metrics that satisfy weaker differential conditions. Let $(M,J,g)$ be a Hermitian manifold of complex dimension $n$ and denote by $\omega$ the fundamental form of $g$. The metric $g$ is called {\em balanced} if $d \omega^{n-1}=0$ \cite{M1982}, {\em strong K\"ahler with torsion} (SKT) if $\del \delbar \omega = 0$ \cite{B1989}, {\em astheno-K\"ahler} if $\del \delbar \omega^{n-2}=0$ \cite{JY1993}, {\em Gauduchon} if $\del \delbar \omega^{n-1}=0$ \cite{G1977}, {\em Strongly Gauduchon} if $\del \omega^{n-1}$ is $\delbar$-exact \cite{P2013}. 

Another class of geometric structures arises in the non-K\"ahler setting, namely, the {\em $p$-K\"ahler}, {\em $p$-symplectic}, and {\em $p$-pluriclosed structures}. These structures were first introduced in \cite{AA1987, A2011, A2017}. A real {\em transverse} $(p,p)$-form is called $p$-K\"ahler if it is $d$-closed, $p$-pluriclosed if it is $\del \delbar$-closed (see \cite{A2017, AA1987}, Definition \ref{Definitions of p-structures}.) A real $d$-closed $2p$-form such that its $(p,p)$-component is transverse is called $p$-symplectic structure (see \cite{AA1987}, Definition \ref{Definitions of p-structures}). We recall that the notion of {\em transversality} was first introduced in \cite{HK1974}, and independently by D. Sullivan (see \cite{S1976}) in the context of cone structures, defined as a continuous field of cones of $p$-vectors on a manifold. 

For $p=1$ and $p=n-1$, $p$-K\"ahler (resp. $p$-symplectic, $p$-pluriclosed) structures coincide with special Hermitian structures. Specifically, $1$-K\"ahler structures coincide with K\"ahler metrics, $1$-symplectic structures with {\em Hermitian-symplectic structures} \cite{ST2010} and $1$-pluriclosed structures with SKT metrics. Similarly, $(n-1)$-K\"ahler structures coincide with balanced metrics, $(n-1)$-symplectic structures coincide with strongly Gauduchon metrics and $(n-1)$-pluriclosed structures coincide with Gauduchon metrics. 

Some examples of compact complex manifolds admitting special structures are given by {\em solvmanifolds} and {\em nilmanifolds}. By a solvmanifold (resp. nilmanifold), we mean a compact quotient of a connected, simply connected solvable (resp. nilpotent) Lie group by a lattice, equipped with a left-invariant complex structure.

In contrast to the K\"ahler case, the non-K\"ahler setting lacks many of the cohomological and stability properties typically associated with K\"ahler manifolds. There are examples of balanced metrics whose $(n-1)$-power of its fundamental form defines the zero class in the de Rham cohomology (see \cite{OUV2017, PT2020}). Moreover, in general, Dolbeault, Bott–Chern, and Aeppli cohomologies are not isomorphic to each other. Special Hermitian metrics are also not stable under small deformations (see \cite{AB1990}, \cite{FT2009}). The lack of stability under deformations extends to p-K\"ahler structures and $p$-pluriclosed structures as well.

In \cite{W2006, AU2017, RWZ2019, RWZ2022} sufficient conditions have been identified for the stability of balanced metrics and $p$-K\"ahler structures under small deformations. These conditions are often tied to the validity of the $\del \delbar$-lemma or to weaker variants thereof. 

It is important to note that for $1 < p < n-1$, $p$-K\"ahler structures do not arise from Hermitian metrics. A $p$-K\"ahler structure cannot be the $p$-th power of the fundamental form of a Hermitian metric unless the metric is K\"ahler. Hence, in the context of non-K\"ahler geometry, $p$-K\"ahler structures have no metric meaning, unless for extremal cases. However, in some cases, the $(p,p)$-component of a $p$-symplectic structure or a $p$-pluriclosed structure may coincide with the $p$-th power of the fundamental form of a Hermitian metric, although this is not true in general.

Nevertheless, the Alessandrini–Bassanelli conjecture (see \cite[pp. 299]{AB1992}) states that on a compact complex manifold, the existence of a $p$-K\"ahler structure implies the existence of a $(p+1)$-K\"ahler structure. A direct consequence of this conjecture would be that $p$-K\"ahler geometry is a special case of balanced geometry.

\medskip
One of the central themes of this paper is the Alessandrini-Bassanelli conjecture, which we address in the context of {\em holomorphically parallelizable nilmanifolds} and of nilmanifolds equipped with {\em nilpotent complex structures}. In \cite[pp. 299]{AB1992}, the authors state that the conjecture holds true for holomorphically parallelizable nilmanifolds, although no proof is provided. In this paper, we provide a proof of this statement under the hyphotesis that the complex structure equations satisfying \eqref{J-nilpotent for holomorphically parallizable} have Gaussian rational coefficients (see Theorem \ref{Alessandrini-Bassanelli conjecture}). 
\begin{theorem}
    Let $(M,J)$ be a holomorphically parallelizable nilmanifold of complex dimension $n$. Suppose that $(M,J)$ admits a basis of $(1,0)$-forms $\{\varphi^{1}, \dots, \varphi^{n}\}$ such that the complex structure equations satisfy 
    \begin{equation*}
        d \varphi^{j} = \sum_{i < k < j} A_{jik} \varphi^{i} \wedge \varphi^{k}, \quad j=1, \dots,n,
    \end{equation*}
    where $A_{jik} \in \Q[i]$. If $(M,J)$ admits a $p$-K\"ahler structure, then it admits a $(p+1)$-K\"ahler structure.
\end{theorem}

\medskip

Furthermore, $(n-2)$-K\"ahler nilmanifolds equipped with nilpotent complex structures are balanced (see Theorem \ref{Alessandrini-Bassanelli conjecture per n-2}).
\begin{theorem}
    Let $(M,J)$ be a nilmanifold of complex dimension $n$ equipped with a nilpotent complex structure $J$. If $(M,J)$ admits a $(n-2)$-K\"ahler structure, then $(M,J)$ is balanced. 
\end{theorem}

\medskip

One of the goals of this paper is to study $p$-K\"ahler structures on differentiable families of compact complex manifolds. In particular, we provide necessary cohomological conditions on the central fiber of such a family for the existence of a smooth curve of $p$-K\"ahler structures (see Theorem \ref{Necessary condition for the existence of curve of 3-K}). 
\begin{theorem}
    Let $(M,J,\Omega)$ be a compact $p$-K\"ahler manifold of complex dimension $n$, for $1 < p < n-1$, and let $\{M_{t}\}_{t \in I}$ be a differentiable family of compact complex manifolds, where $I \doteq (-\epsilon, \epsilon), \epsilon >0$ and $M_{0} \doteq M$. Suppose that $\{M_{t}\}_{t \in I}$ is parametrized by the $(0,1)$-vector form $\varphi(t) \in \mathcal{A}^{0,1}(T^{1,0} M)$.

    Let $\{\Omega_{t}\}_{t \in I}$ be a smooth family of real transverse $(p,p)$-forms along $\{M_{t}\}_{t \in I}$ as in \eqref{smooth family of (p,p)-forms}. If every $\Omega_{t}$ is a $p$-K\"ahler structure, then 
    \begin{equation*}
        \del \circ \iota_{\varphi^{'}(0)}(\Omega) = - \delbar \, \Omega^{'}(0).
    \end{equation*}
\end{theorem}

\medskip
The techniques used to prove this result are the same as those used in \cite{PS2021, S2022, S2023} for the study of balanced, SKT and astheno-K\"ahler metrics. We remark that case $p \doteq n-1$ is studied in \cite{S2022}.

Moreover, a direct consequence is the following obstruction involving the $(p,p+1)$-Dolbeault group of the central fiber (see Corollary \ref{Obstruction in dolbeault cohomology}).
\begin{corollary}
    Let $(M,J,\Omega)$ be a compact $p$-K\"ahler manifold of complex dimension $n$, for $1 < p < n-1$, and let $\{M_{t}\}_{t \in I}$ be a differentiable family of compact complex manifolds, where $I \doteq (-\epsilon, \epsilon), \epsilon >0$ and $M_{0} \doteq M$. Suppose that $\{M_{t}\}_{t \in I}$ is parametrized by the $(0,1)$-vector form $\varphi(t) \in \mathcal{A}^{0,1}(T^{1,0} M)$. If there exists a smooth family of $p$-K\"ahler structures $\{\Omega_{t}\}_{t \in I}$ such that $\Omega_{0} \doteq \Omega$, then 
    \begin{equation*}
        \Big[ \del \circ \iota_{\varphi^{'}(0)} (\Omega) \Big]_{H^{p,p+1}_{\delbar}(M)} = 0.
    \end{equation*}
\end{corollary}
We remark that the obstruction results above are related to the existence of smooth families of $p$-K\"ahler structures along differentiable families of compact complex manifolds, once a fixed $p$-K\"ahler structure on the central fiber is given.
\medskip

Finally, we address the problem of studying the cohomology classes of $p$-K\"ahler, $p$-symplectic and $p$-pluriclosed structures. In general, a $p$-K\"ahler structure $\Omega$ defines a cohomology class $[\Omega]_{\#}$, where $\# = \text{BC}, \del, \delbar, \text{A}, \text{dR}$. A natural question is to determine necessary conditions for the vanishing of the cohomology class of $\Omega$ (see Theorem \ref{exactness of p-Kahler form}).
\begin{theorem}
    Let $(M,J,\Omega)$ be a compact $p$-K\"ahler manifold of complex dimension $n$, for \\
    $1 < p \leq n-1$. Let $k \doteq n-p$. If 
    \begin{equation*}
        [\Omega]_{\#} = 0, \quad \text{for $\# = dR, A, BC, \partial, \overline{\partial}$},  
    \end{equation*}
    then, there are no non-vanishing, simple, $(k,0)$-forms $\xi$ such that $\del \xi = \delbar \xi = 0$.
\end{theorem}

Similar results are proved for $p$-symplectic and $p$-pluriclosed structures (see Theorems \ref{exactness of p-sympl}, \ref{Necessary condition for the vanishing of p-pluriclosed forms}). 
\begin{theorem}
   Let $(M,J,\Psi)$ be a compact $p$-symplectic manifold of complex dimension $n$, for $1 \leq p \leq n-1$. Let $k \doteq n-p$. If 
    \begin{equation*}
        [\Psi]_{dR} =0 ,  
    \end{equation*}
    then, there are no non-vanishing, simple, $(k,0)$-forms $\xi$ such that $\del \xi = \delbar \xi = 0$.
\end{theorem}
\begin{theorem}
    Let $(M,J,\Omega)$ be a compact $p$-pluriclosed manifold of complex dimension $n$, for $1 \leq p \leq n-1$. Let $k \doteq n-p$. If 
    \begin{equation*}
        [\Omega]_{A} =0 ,  
    \end{equation*}
    then, there are no non-vanishing, simple, $(k,0)$-forms $\xi$ such that $\del \xi = \delbar \xi = 0$.
\end{theorem}

\medskip

We note that in \cite{PT2020}, the authors study necessary cohomological conditions for the vanishing of the Aeppli cohomology class of the $p$-th power of the fundamental form of a Hermitian metric. 

\medskip

In addition to the theoretical results, we construct several explicit examples of $p$-K\"ahler and $p$-symplectic compact complex manifolds. We recall that in \cite{AB1991, A1999}, some examples of compact $p$-K\"ahler manifolds are provided. Specifically, in \cite{AB1991} it is shown that the compact holomorphically parallelizable nilmanifolds $\eta \beta_{2n + 1}$, which generalize the Iwasawa manifold, admit $p$-K\"ahler structures for $n+1 \leq p \leq 2n+1$. On the other hand, in \cite{A1999}, the author provides some examples of Moishezon manifolds that admit $p$-K\"ahler structures. 

Here, we construct a $3$-K\"ahler structure on $\eta \beta_{5}$ that is different from the one constructed by L. Alessandrini and G. Bassanelli in \cite{AB1991} (Example \ref{etabeta_{5}}). We note that this construction is made possible by the results established in \cite{FM2025}, where the authors provide methods to construct transverse forms. Moreover, we provide a family of $3$-K\"ahler nilmanifolds of complex dimension $5$ equipped with nilpotent complex structures (Example \ref{3-kahler nilmanifolds}). We also exhibit a family of $3$-symplectic nilmanifolds of complex dimension $5$ equipped with nilpotent complex structures (Example \ref{3-symplectic nilmanifolds}). Finally, we construct a smooth family of $3$-K\"ahler structures along a differentiable family of compact complex manifolds whose central fiber is $\eta \beta_{5}$ (Example \ref{examples of deformations}).

\medskip

The paper is organized as follows. In Section \ref{Section 2} we establish the notation that is adopted throughout the paper and recall the definitions of $p$-K\"ahler, $p$-symplectic and $p$-pluriclosed structures. In Section \ref{Section 3} we recall some basic facts on the theory of deformations of complex manifolds, in particular, we recall the notation adopted in \cite{RZ2015, RZ2018}.

In Section \ref{Section 4}, we address the conjecture of L. Alessandrini and G. Bassanelli on holomorphically parallelizable nilmanifolds and nilmanifolds equipped with nilpotent complex structures. Moreover, in this section we construct examples of compact $p$-K\"ahler manifolds and compact $p$-symplectic manifolds. 

Section \ref{Section 5} is devoted to the study of necessary cohomological conditions for the existence of curves of $p$-K\"ahler structures along differentiable families of compact complex manifolds starting from a fixed $p$-K\"ahler structure. We also study two different curves of $p$-K\"ahler structures starting from different $p$-K\"ahler structures on $\eta \beta_{5}$.

Finally, in Section \ref{Section 6} we study the cohomology classes associated with $p$-K\"ahler, $p$-symplectic and $p$-pluriclosed structures. As an application, we show that on nilmanifolds equipped with left-invariant complex structure, the cohomology classes of $p$-K\"ahler, $p$-symplectic and $p$-pluriclosed structures are always non-vanishing.

\vskip.3truecm
{\em \underline{Acknowledgments:}} The author would like to sincerely thank Adriano Tomassini for his support, encouragement, and many valuable discussions and suggestions. Special thanks also go to Asia Mainenti for her careful reading of the paper and many helpful discussions. The author is also grateful to Anna Fino, Sheng Rao and Lorenzo Sillari for their helpful comments and remarks. Many thanks are also due to Riccardo Piovani for suggesting some possible further developments and to Marco Andreatta for bringing my attention to the paper \cite{A1999}.

\section{Preliminaries on p-K\"ahler, p-symplectic and p-pluriclosed structures}\label{Section 2}
Given a complex manifold $(M,J)$ of complex dimension $n$, we denote by $T_{\C}M \doteq TM \otimes \C$ its complexified tangent bundle. The complexified tangent bundle $T_{\C}M$ can be decomposed as $T_{\C}M = T^{1,0}M \oplus T^{0,1}M$, where $T^{1,0}M \doteq \{X \in T_{\C}M \, | \, JX =i X\}$ and $T^{0,1}M \doteq \{X \in T_{\C}M \, | \, JX =-i X\}$. The bundle $ \Lambda^{r}_{\C}M$ of complex $r$-forms splits, consequently, as 
\begin{equation*}
    \Lambda^{r}_{\C}M = \oplus_{p+q=r} \Lambda^{p,q}M,
\end{equation*} 
where $\Lambda^{p,q}M \doteq \Lambda^{p}(T^{1,0}M)^{\ast} \otimes \Lambda^{q}(T^{0,1}M)^{\ast}$.

A form $\alpha \in \Lambda^{p,p} (M)$ (resp. $\alpha \in \Lambda^{r}_{\C}(M)$) is called real if $\alpha = \overline{\alpha}$. We adopt the following notation 
\begin{equation*}
    \Lambda^{r}_{\R} M \doteq \{ \alpha \in \Lambda^{r}_{\C}M \, | \, \text{$\alpha$ is real}\}, \quad \Lambda^{p,p}_{\R} M \doteq \{\alpha \in \Lambda^{p,p}M \, | \, \text{$\alpha$ is real} \}.
\end{equation*}
We denote by $\mathcal{A}^{p,q}(M)$ and $\mathcal{A}^{r}(M,\C)$ the spaces of smooth sections of the bundles $\Lambda^{p,q}(M)$ and $\Lambda^{r}_{\C}(M)$, respectively.

\medskip

Throughout the following discussion, we will work pointwise. 

Let $T_{x}^{\ast}M$ be the dual of the tangent space of $M$ at the point $x$. Fix a basis $\{\varphi^{j}\}_{j=1}^{n}$ for the space $\Lambda^{1,0}(T_{x}^{\ast}M \otimes \C)$. Then, a basis for $\Lambda^{p,q}(T_{x}^{\ast}M \otimes \C)$ is given by 
\begin{equation*}
    \{\varphi^{i_{1}} \wedge \dots \wedge \varphi^{i_{p}} \wedge \varphi^{\overline{j_{1}}} \wedge \dots \wedge \varphi^{\overline{j_{q}}} \; | \; 1 \leq i_{1} < \dots < i_{p} \leq n, 1 \leq j_{1} < \dots < j_{p} \leq n \}.
\end{equation*}
A volume form is given by
\begin{equation*}
    \text{Vol} \doteq (\frac{i}{2} \varphi^{1} \wedge \overline{\varphi}^{1}) \wedge \dots \wedge (\frac{i}{2} \varphi^{n} \wedge \overline{\varphi}^{n}) = \sigma_{n} \varphi^{1}\wedge \dots \wedge \varphi^{n} \wedge \overline{\varphi}^{1} \wedge \dots \wedge \overline{\varphi}^{n},
\end{equation*}
where $\sigma_{p} \doteq i^{p^{2}}2^{-p}$.

An element $\tau \in \Lambda^{n,n}_{\R}(T^{\ast}_{x}M \otimes \C)$, is said to be {\em positive} (resp.{\em  strictly positive}) if $\tau = a \text{Vol}$, where $a \geq 0$ ($a > 0$).

\begin{definition}
     A form $\psi \in \Lambda^{p,0}(T^{\ast}_{x}M \otimes \C)$ is said to be {\em simple} if $ \psi = \psi^{1} \wedge \dots \wedge \psi^{p}$, where $\psi^{i} \in \Lambda^{1,0}(T^{\ast}_{x}M \otimes \C)$, for $i=1, \dots, p$.
\end{definition}

Let $\Omega \in  \Lambda^{p,p}_{\R}(T^{\ast}_{x}M \otimes \C)$. We say that $\Omega$ is {\em transverse} ($\Omega > 0$), if 
\begin{equation*}
    \sigma_{n-p} \Omega \wedge \beta \wedge \overline{\beta} 
\end{equation*}
is strictly positive, $\forall \beta \in \Lambda^{n-p,0}(T^{\ast}_{x}M \otimes \C)$ simple, $\beta \neq 0$.

\begin{definition}\label{Definitions of p-structures}
    Let $\Omega$ be a real $(p,p)$-form on $M$ such that $\Omega_{x} \in \Lambda^{p,p}_{\R}(T_{x}^{\ast}M \otimes \C)$ is transverse $\forall x \in M$, and let $\Psi$ be a real $2p$-form on $M$ such that $\Psi^{p,p} \doteq \Omega$. 
    \begin{enumerate}
        \item $\Omega$ is called a $p$-{\em K\"ahler} structure if it is $d$-closed (the triple $(M,J,\Omega)$ is called a {\em $p$-K\"ahler manifold});
        \item $\Omega$ is called a $p$-{\em pluriclosed} structure if it is $\partial \overline{\partial}$-closed (the triple $(M,J,\Omega)$ is called a {\em $p$-pluriclosed manifold});
        \item $\Psi$ is called a $p$-{\em symplectic} structure if it is $d$-closed (the triple $(M,J,\Psi)$ is called a {\em $p$-symplectic manifold}).
    \end{enumerate}
\end{definition}

We recall that a complex manifold $(M,J)$ of complex dimension $n$ is called {\em holomorphically parallelizable} if $T^{1,0}M$ is trivial as a holomorphic bundle, namely, there exist $n$ holomorphic $1$-forms on $(M,J)$ which are linearly independent at every point of $M$. In \cite{W1954}, the author shows that a compact complex manifold is holomorphically parallelizable if and only if it is the compact quotient of a connected, simply-connected, complex Lie group by a lattice.

Let $(\Gamma \backslash G,J)$ be a nilmanifold and denote by $\mathfrak{g}$ the Lie algebra of $G$. The complex structure $J$ is called {\em nilpotent} if the ascending series $\{\mathfrak{g}_{i}^{J}\}$ defined by
\begin{equation*}
    \mathfrak{g}^{J}_{0} \doteq 0, \quad \mathfrak{g}^{J}_{i} \doteq \{X \in \mathfrak{g} \enskip | \enskip [X, \mathfrak{g}]\subseteq \mathfrak{g}^{J}_{i-1}, \quad [J X, \mathfrak{g}] \subseteq \mathfrak{g}^{J}_{i-1} \}
\end{equation*}
satisfies $\mathfrak{g}^{J}_{k}= \mathfrak{g}$ for some $k > 0$. An equivalent condition (see \cite{CFGU1997}, \cite{CFGU2000}) is that there exists a coframe of left-invariant $(1,0)$-forms $\{\varphi^{i}\}_{i=1, \dots n}$  such that 
\begin{equation}\label{J-nilpotent}
    d \varphi^{j} = \sum_{i<k<j} A_{jik} \varphi^{i k} + \sum_{i,k <j} B_{jik} \varphi^{i \overline{k}}, \quad j=1,\dots n,
\end{equation}
where $A_{jik}, B_{jik} \in \C$. We recall that holomorphically paralellizable nilmanifolds are such that the coefficients $B_{jik}$ in \eqref{J-nilpotent} vanish for every $j=1, \dots, n$ (see \cite{CFGU1997}, \cite{CFGU2000}), i.e., 
\begin{equation}\label{J-nilpotent for holomorphically parallizable}
    d \varphi^{j} = \sum_{i<k<j} A_{jik} \varphi^{i k}, \quad j=1,\dots n,
\end{equation}
where $A_{jik} \in \C$.

\section{Preliminaries on deformations theory}\label{Section 3}
In this section we recall the fundamental facts on deformations theory and establish the notation that will be adopted throughout the paper. 

\begin{definition}
    Let $B$ be a connected differentiable manifold of real dimension $m$ and let $\{M_{t}\}_{t \in B}$ be a set of compact complex manifolds. We say that $\{M_{t}\}_{t \in B}$ is a {\em differentiable family of compact complex manifolds} if there exists a differentiable manifold $\mathfrak{M}$ and a differentiable, proper map $\pi : \mathfrak{M} \to B$ such that 
    \begin{enumerate}
        \item $\pi^{-1}(t) = M_{t}$ as complex manifold, $\forall t \in B$;
        \item the rank of the Jacobian of $\pi$ is constant and equal to the dimension of $B$, for every $p \in \mathfrak{M}$.
    \end{enumerate}
    If $M_{t_{0}} = M$, we say that $M_{t}$ is a deformation of $M$,  $\forall t \in B$.
\end{definition}

\begin{remark}
    We recall that, by a classical result of C. Ehresmann, if $\{M_{t}\}_{t \in B}$ is a differentiable family of compact complex manifolds, then $M_{t}$ is diffeomorphic to $M_{\widetilde{t}}$, $\forall t, \widetilde{t} \in B$ (see \cite{E1947}, \cite[Proposition 6.2.2]{H2005}). 
\end{remark}
For the sake of simplicity, we assume that $B \doteq \{(t_{1}, \dots, t_{m}) \in \R^{m} \, |  \, |t_{j}| < \epsilon, j =1, \dots m\}\subseteq \R^{m}$, for $\epsilon > 0$ and that $t_{0} = 0$. The compact complex manifold $(M_{0},J_{0})$ is referred to as the {\em central fiber} of the differentiable family of compact complex manifolds. Moreover, if the real dimension of $B$ is $1$, i.e., $B =I =(-\epsilon, \epsilon)$, for $\epsilon>0$, then we refer to the differentiable family $\{M_{t}\}_{t \in I}$ as a {\em curve of deformations} of $(M_{0},J_{0})$.

\medskip

Let $(M,J)$ be a compact complex manifold of complex dimension $n$ and let $\{M_{t}\}_{t \in B}$ be a differentiable family of compact complex manifolds such that $(M_{0},J_{0}) \doteq (M,J)$. In \cite[Chapter 4, Section 4.1 (b)]{K2005} and \cite{KM2006}, is shown that the complex structure on each $M_{t}$, for $t \in B$, is parametrized by a $(0,1)$-vector form $\varphi(t)$ on the central fiber, i.e., $\varphi(t) \in \mathcal{A}^{0,1}(T^{1,0}(M))$. In particular, $\varphi(0) = 0$ and $\varphi(t)$ satisfies the Maurer-Cartan equation, i.e.,
\begin{equation*}
    \delbar \varphi(t) - \frac{1}{2}[\varphi(t),\varphi(t)] = 0,
\end{equation*}
where $[\cdot, \cdot]$ is a suitable bracket. Sometimes, we omit the dependence on $t$ of the $(0,1)$-vector form $\varphi(t)$. The local expression of $\varphi(t)$ (see \cite[Chapter 4]{KM2006}) is given by 
\begin{equation*}
    \varphi(t) = \sum_{\lambda = 1}^{n} \varphi^{\lambda} \otimes \frac{\del}{\del z^{\lambda}},
\end{equation*}
where $\varphi^{\lambda}$ is a global $(0,1)$-form on $(M,J)$ and $(z_{1}, \dots, z_{n})$ are local holomorphic coordinates on $(M,J)$. 

We also recall that a differentiable function $f$ defined on any open subset of $M$ is holomorphic with respect to the complex structure of $M_{t}$ if and only if 
\begin{equation*}
    \big(\delbar -  \sum_{\lambda=1}^{n} \varphi^{\lambda} \otimes \frac{\del}{\del z^{\lambda}} \big) f  = 0,
\end{equation*}
(see \cite[Chapter 4, Proposition 1.2]{KM2006}).

\medskip

With the aid of the $(0,1)$-vector form $\varphi(t)$, we are able to link the $(p,q)$-form on the central fiber to the $(p,q)$-form on any fiber. In the following, we will adopt the notation of \cite{RZ2015, RZ2018}. For the sake of completeness, we recall the notation and some results. 

In general, if $\phi$ is a $(0,s)$-vector form, i.e., $\phi = \eta \otimes Z \in \mathcal{A}^{0,s}(T^{1,0}(M))$, where $\eta \in \mathcal{A}^{0,s}(M)$, $Z \in \Gamma(M,T^{1,0}(M))$, we define a {\em contraction operator} that acts in the following way 
\begin{equation*}
    \iota_{\phi} : \mathcal{A}^{p,q}(M) \to \mathcal{A}^{p-1, q+s}(M), \quad \sigma \mapsto \iota_{\phi}(\sigma) = \eta \wedge (\iota_{Z} \sigma),
\end{equation*}
where $\iota_{Z}$ is the usual contraction operator. Sometimes, we use the notation $\phi \lrcorner \, \sigma$ to indicate the contraction expressed above.

Let $\varphi(t)$ be the $(0,1)$-vector form that parametrizes a differentiable family of compact complex manifolds $\{M_{t}\}_{t \in B}$ such that $(M,J)$ is the central fiber. Since $\varphi(t) \in \mathcal{A}^{0,1}(T^{1,0}(M))$, we are able to define the following operators
\begin{equation}\label{exponential of phi}
    e^{\iota_{\varphi(t)}} \doteq \sum_{k=0}^{\infty} \frac{1}{k!} \iota^{k}_{\varphi(t)}, \quad e^{\iota_{\overline{\varphi(t)}}} \doteq \sum_{k=0}^{\infty} \frac{1}{k!} \iota^{k}_{\overline{\varphi(t)}},
\end{equation}
where $\iota^{k}_{\varphi(t)}$ denotes the contraction by $\varphi(t)$ applied $k$-times, and $\overline{\varphi(t)} \in \mathcal{A}^{1,0}(T^{0,1}(M))$ is the conjugate of $\varphi(t)$. We notice that the sums in the equations above are finite due to the finiteness dimension of $M$.

We now define the following map, first introduced in \cite{RZ2015}, which extends the maps in \eqref{exponential of phi}:
\begin{equation*}
    e^{\iota_{\varphi(t)}|\iota_{\overline{\varphi(t)}}} : \mathcal{A}^{p,q}(M) \to \mathcal{A}^{p,q}(M_{t}).
\end{equation*}
This map acts as follows. Given $\sigma \in \mathcal{A}^{p,q}(M)$, locally written as $\sigma = \sigma_{i_{1} \dots i_{p} j_{1} \dots j_{q}} dz^{i_{1}} \wedge \dots \wedge dz^{i_{p}} \wedge dz^{\overline{j_{1}}} \wedge \dots \wedge dz^{\overline{j_{q}}}$, then 
\begin{equation*}
    e^{\iota_{\varphi(t)}|\iota_{\overline{\varphi}}} (\sigma) = \sigma_{i_{1} \dots i_{p} j_{1} \dots j_{q}}  e^{\iota_{\varphi(t)}}(dz^{i_{1}} \wedge \dots \wedge dz^{i_{p}}) \wedge (e^{\iota_{\overline{\varphi(t)}}} (dz^{\overline{j_{1}}} \wedge \dots \wedge dz^{\overline{j_{q}}}) ).
\end{equation*}
It turns out (see \cite[Lemma 2.9]{RZ2018}) that $e^{\iota_{\varphi(t)}|\iota_{\overline{\varphi(t)}}}$ is a linear isomorphism for $t$ arbitrarily small.

We also introduce the following notation: 
\begin{equation*}
    \varphi(t) \Finv \, \sigma \doteq \sigma_{i_{1} \dots i_{p} j_{1} \dots j_{q}} \iota_{\varphi(t)} (dz^{1}) \wedge \dots \wedge \iota_{\varphi(t)} ( dz^{i_{p}}) \wedge \iota_{\varphi(t)} (dz^{\overline{j_{1}}})\wedge \dots \wedge \iota_{\varphi(t)}( dz^{\overline{j_{q}}}).
\end{equation*}
This contraction is called {\em simultaneous contraction}, and it allows us to rewrite the operator $e^{\iota_{\varphi(t)}|\iota_{\overline{\varphi(t)}}}$ as 
\begin{equation*}
    e^{\iota_{\varphi(t)}|\iota_{\overline{\varphi(t)}}} = \big(\text{I} + \varphi(t) + \overline{\varphi}(t) \big) \Finv.
\end{equation*}
In \cite[Proposition 2.7, Proposition 2.13]{RZ2018} is shown that the differential operators $\del_{t}, \delbar_{t}$ on each $M_{t}$, for $t \in B$, can be expressed in terms of formulas involving $\varphi(t)$. Indeed, if $f$ is a differentiable function on $M$, with values in $\C$, we have that 
\begin{equation*}
    \begin{split}
        & \del_{t} f = e^{\iota_{\varphi}} \Big(\big(I - \varphi \overline{\varphi} \big)^{-1} \lrcorner \, \big(\del - \overline{\varphi} \lrcorner \, \delbar \big) f \Big), \\
        & \delbar_{t} f =  e^{\iota_{\overline{\varphi}}} \Big(\big(I - \overline{\varphi} \varphi  \big)^{-1} \lrcorner \, \big(\delbar - \varphi \lrcorner \, \del \big) f \Big),
    \end{split}
\end{equation*}
where $\varphi \overline{\varphi} \doteq \overline{\varphi} \lrcorner \, \varphi$, $\overline{\varphi} \varphi \doteq \varphi \lrcorner \, \overline{\varphi}$ and $\varphi \doteq \varphi(t)$. Furthermore, the actions of $\del_{t}, \delbar_{t}$ on a $(p,q)$-form on $M_{t}$ written as $e^{\iota_{\varphi}|\iota_{\overline{\varphi}}} \alpha$ for $\alpha \in \mathcal{A}^{p,q}(M)$ are the following  
\begin{equation}\label{del operator on M_t}
    \del_{t} \big(e^{\iota_{\varphi}|\iota_{\overline{\varphi}}} \alpha \big) = e^{\iota_{\varphi}|\iota_{\overline{\varphi}}} \Big(\big(I - \varphi \overline{\varphi} \big)^{-1} \Finv \, \big([\delbar,\iota_{\overline{\varphi}}] + \del \big) (I - \varphi \overline{\varphi} \big) \Finv \, \alpha \Big),
\end{equation}
\begin{equation}\label{delbar operator on M_t}
    \delbar_{t} \big(e^{\iota_{\varphi}|\iota_{\overline{\varphi}}} \alpha \big) = e^{\iota_{\varphi}|\iota_{\overline{\varphi}}} \Big(\big(I - \overline{\varphi} \varphi \big)^{-1} \Finv \, \big([\del,\iota_{\varphi}] + \delbar \big) (I - \overline{\varphi} \varphi \big) \Finv \, \alpha \Big).
\end{equation}

\section{p-K\"ahler and p-symplectic structures on nilmanifolds}\label{Section 4}
In this section, we address the conjecture of L. Alessandrini and G. Bassanelli (see \cite[pp. 299]{AB1992}) on holomorphically parallelizable nilmanifolds and $(n-2)$-K\"ahler nilmanifolds equipped with nilpotent complex structures. In the end of the section, we study some examples. In particular, we construct a new $3$-K\"ahler structure on $\eta \beta_{5}$ and we provide a family of $3$-K\"ahler nilmanifolds of complex dimension $5$. In the last example, we construct a family of nilmanifolds of complex dimension $5$ that admit astheno-K\"ahler metrics and $3$-symplectic structures.

\begin{theorem}\label{Alessandrini-Bassanelli conjecture}
    Let $(M,J)$ be a holomorphically parallelizable nilmanifold of complex dimension $n$. Suppose that $(M,J)$ admits a basis of $(1,0)$-forms $\{\varphi^{1}, \dots, \varphi^{n}\}$ such that the complex structure equations satisfy 
    \begin{equation*}
        d \varphi^{j} = \sum_{i < k < j }A_{jik} \varphi^{i} \wedge \varphi^{k}, \quad j=1, \dots,n,
    \end{equation*}
    where $A_{jik} \in \Q[i]$. If $(M,J)$ admits a $p$-K\"ahler structure, then it admits a $(p+1)$-K\"ahler structure.
\end{theorem}
\begin{proof}
    Recall that every left-invariant metric on $(M,J)$ is balanced (see \cite{AG1986}); thus, every holomorphically parallelizable nilmanifold admitting a $(n-2)$-K\"ahler structure also admits a $(n-1)$-K\"ahler structure. Therefore, the statement holds for complex dimensions $n \leq 4$. We now assume that $p < n - 2$ and $k \doteq n-p$.

    \medskip
    
    We prove the theorem by induction on the complex dimension $n$ of the manifold. As base case we consider $n=5$. Since the complex structure of a holomorphically parallelizable nilmanifold is nilpotent, in complex dimension $5$, the thesis follows by \cite[Theorem 3.8]{FM2024}.

    \medskip
    
    Assuming the statement holds for complex dimension $m \doteq n-1$, we aim to prove it for $n \geq 6$. Our approach is as follows: first, we show that the existence of a $p$-K\"ahler structure imposes necessary conditions on the complex structure equations of $(M,J)$; then, we prove that if $(M,J)$ does not admit a $(p+1)$-K\"ahler structure, it cannot admit a $p$-K\"ahler structure either.

    \medskip
    
    Let us denote by $G$ the simply connected, connected, complex nilpotent Lie group whose compact quotient by a discrete subgroup is $M$ (see \cite{W1954}). By hypothesis, the complex structure equations of the Lie algebra $\mathfrak{g}$ of $G$, endowed with the complex structure induced by that of $(M,J)$, are 
    \begin{equation*}
        \begin{cases}
            d \varphi^{1} = d \varphi^{2} = 0, \\
            d\varphi^{r} = \sum_{u < v < r}A_{ruv} \varphi^{u} \wedge \varphi^{v}, \quad r=3, \dots,n,
        \end{cases}
    \end{equation*}
    where $A_{ruv} \in \Q[i]$. We consider $\mathfrak{b}_{\C}^{\ast} = \langle \varphi^{n}, \varphi^{\overline{n}} \rangle$ and  $\mathfrak{k}_{\C}^{\ast}=\langle\varphi^{1}, \dots, \varphi^{n-1}, \varphi^{\overline{1}}, \dots, \varphi^{\overline{n-1}} \rangle$. The complex Lie algebra $\mathfrak{k}$, whose complex dimension is $m = n-1$, admits a $(p-1)$-K\"ahler structure (see \cite[Proof of Proposition 3.3]{FM2024}). 
    
    Given that $k=n-p$, for greater clarity in what follows, we find it convenient to state that $\mathfrak{k}$ admits a $(m-k)$-K\"ahler structure. 
    
    The complex structure equations of $\mathfrak{k}$ are 
    \begin{equation*}
        \begin{cases}
            d \varphi^{1} = d \varphi^{2} = 0, \\
            d\varphi^{r} = \sum_{u < v < r}A_{ruv} \varphi^{u} \wedge \varphi^{v}, \quad r=3, \dots,n-1,
        \end{cases}
     \end{equation*}
     where $A_{ruv} \in \Q[i]$. Hence, due to Malcev's Theorem \cite{M1951}, we are sure that the simply connected, connected, complex Lie group $K$, whose Lie algebra is $\mathfrak{k}$ admits a lattice $\Gamma$, hence $(M_{\mathfrak{k}} \doteq\Gamma \backslash K,J_{\mathfrak{k}})$ is a holomorphically parallelizable nilmanifold. Since the complex dimension of $M_{\mathfrak{k}}$ is $m$, we can use the induction hyphotesis to conclude that $(M_{\mathfrak{k}},J_{\mathfrak{k}})$ admits $(m-l)$-K\"ahler structures for $l= 0,1, \dots, k$.

    \medskip
    
    Note that, by \cite[Theorem 2.3]{ST2022}, since $(M_{\mathfrak{k}},J_{\mathfrak{k}})$ admits $(m-l)$-K\"ahler structures for $l= 0,1, \dots, k$ and $p < n - 2$, we have that, up to re-ordering,
     \begin{equation}
         d \varphi^{j} = 0, \quad j = 1, \dots,k + 1.
     \end{equation}
     Hence, the complex structure equations of $M$ are 
     \begin{equation*}
        \begin{cases}
            d \varphi^{j}=0, \quad j=1, \dots, k+1, \\
            d\varphi^{r} = \sum_{u < v < r}A_{ruv} \varphi^{u} \wedge \varphi^{v},  \quad r= k+2, \dots, n,
        \end{cases}
     \end{equation*}
     where $A_{ruv} \in \Q[i]$.

    \medskip
    
    Now, suppose that $(M,J)$ does not admit a $(p+1)$-K\"ahler structure. By \cite[Theorem 3.2]{AB1991}, $(M,J)$ does not admit a $(p+1)$-K\"ahler structure if and only if 
    \begin{equation}
        \exists \, \zeta \in \mathcal{A}^{k-1,0}(M), \quad \text{$\zeta$ is simple}, \quad \delbar \zeta = 0, \quad \zeta = \del \alpha, \quad \text{for $\alpha \in \mathcal{A}^{k-2,0}(M)$}.
    \end{equation}
    First of all notice that $\zeta$ is left-invariant, hence it can be written as
    \begin{equation*}
        \zeta = \sum_{|S|=k-1} a_{S}\varphi^{S},  
    \end{equation*}
    where $S$ varies through the multi-indices $(s_1,\dots,s_{k-1})$ of length $k-1$ and $a_S\in \C$.
    Let us now define $\zeta_{1} \doteq - \varphi^{1} \wedge \zeta$, then 
    \begin{equation*}
        \zeta_{1} \in \mathcal{A}^{k,0}(M), \quad \text{$\zeta_{1}$ is simple}, \quad \delbar \zeta_{1} = 0, \quad \zeta_{1} = \del (\varphi^{1} \wedge \alpha).
    \end{equation*}
    The only issue is that $\zeta_{1}$ could be zero. Notice that $\zeta_{1}$ is zero if and only if 
    \begin{equation}\label{first step}
        \zeta = \varphi^{1} \wedge \sum_{|S|=k-2} a^{1}_{S}\varphi^{S},
    \end{equation}
    where $S$ varies through the multi-indices $(s_1,\dots,s_{k-2})$ of length $k-2$ and $a^{1}_{S} \in \C$.
    If equation \eqref{first step} holds, define $\zeta_{2} \doteq - \varphi^{2} \wedge \zeta$, then 
    \begin{equation*}
        \zeta_{2} \in \mathcal{A}^{k,0}(M), \quad \text{$\zeta_{2}$ is simple}, \quad \delbar \zeta_{2} = 0, \quad \zeta_{2} = \del (\varphi^{2} \wedge \alpha).
    \end{equation*}
    As before, $\zeta_{2}$ could be zero. Notice that, $\zeta_{2}$ is zero if and only if 
    \begin{equation}\label{second step}
        \zeta_{2} = \varphi^{12} \wedge \sum_{|S|=k-3} a^{2}_{S}\varphi^{S},
    \end{equation}
    where $S$ varies through the multi-indices $(s_1,\dots,s_{k-3})$ of length $k-3$ and $a^{2}_{S}\in \C$. We iterate this construction. 
    
    The final two steps are the following. Define
    \begin{equation*}
            \zeta_{k-1}  \doteq - \varphi^{k-1} \wedge \zeta = - \varphi^{k-1} \wedge \varphi^{1 \dots k-2} \wedge \Big(\sum_{s=k-1}^{n} a^{k-1}_{s} \varphi^{s} \Big),
    \end{equation*}
    where $a^{k-1}_{s} \in \C$, thus 
    \begin{equation*}
        \zeta_{k-1} \in \mathcal{A}^{k,0}(M), \quad \text{$\zeta_{k-1}$ is simple}, \quad \delbar \zeta_{k-1} = 0, \quad \zeta_{k-1} = \del(\varphi^{k-1} \wedge \alpha).
    \end{equation*}
    Moreover, $\zeta_{k-1} = 0$ if and only if 
     \begin{equation}\label{zetaj = 0} 
         \zeta_{k-1} = a^{k} \, \varphi^{1} \wedge \varphi^{2} \wedge \dots \wedge \varphi^{k-1}.
     \end{equation}
    Thus, define $\zeta_{k} \doteq - \varphi^{k} \wedge \zeta_{k-1}$, then 
    \begin{equation}\label{zeta2 = 0}
        \zeta_{k} \in \mathcal{A}^{k,0}(M), \quad  \text{$\zeta_{k}$ is simple}, \quad \delbar \zeta_{k} = 0,  \quad \zeta_{k} \neq 0,\quad \zeta_{k} = \del(\varphi^{k} \wedge \alpha).
    \end{equation} 
    Finally, by \cite[Theorem 3.2]{AB1991}, there are no $p$-K\"ahler structures on $(M,J)$ and this is absurd.
\end{proof}

In the following, we show that $(n-2)$-K\"ahler nilmanifolds equipped with nilpotent complex structures are balanced.
\begin{theorem}\label{Alessandrini-Bassanelli conjecture per n-2}
    Let $(M,J)$ be a nilmanifold of complex dimension $n$ equipped with a nilpotent complex structure $J$. If $(M,J)$ admits a $(n-2)$-K\"ahler structure, then $(M,J)$ is balanced. 
\end{theorem}
\begin{proof}
    First of all notice that for complex dimension $n \leq 4$ the thesis is always true. The proof for complex dimension $4$ is given in \cite[Theorem 3.7]{FM2024}. Let us suppose that $n \geq 5$ and let us denote by $\Omega$ the $p$-K\"ahler structure on $(M,J)$.

    \medskip
    
    Requiring that $J$ is nilpotent is equivalent to the existence of a basis $\{\varphi^{1}, \dots, \varphi^{n}\}$ of left-invariant $(1,0)$-forms, such that 
    \begin{equation}\label{J-nilpotent nel th per le (n-2)}
        \begin{cases}
            d \varphi^{1}= 0, \\
            d\varphi^{r} = \sum_{u < v < r}A_{ruv} \varphi^{u} \wedge \varphi^{v} + \sum_{u, v < r}B_{ruv} \varphi^{u} \wedge \varphi^{\overline{v}}, \quad r=2, \dots, n,
        \end{cases}
    \end{equation}
    where $A_{ruv}, B_{ruv} \in \C$.
    
   Note that, if $(M,J)$ admits a $(n-2)$-K\"ahler structure, then there exists a left-invariant $(1,0)$-coframe $\{\psi_{1}, \dots, \psi_{n}\}$ such that 
    \begin{equation}\label{Correct coframe}
        \begin{cases}
            d \psi^{1} = d \psi^{2} = d \psi^{3} = 0, \\
            d \psi^{r} = \sum_{u < v < r}D_{ruv} \psi^{u} \wedge \psi^{v} + \sum_{u, v < r}E_{ruv} \psi^{u} \wedge \psi^{\overline{v}}, \quad r=4, \dots, n,
        \end{cases}
    \end{equation}
    where $D_{ruv}, E_{ruv} \in \C$. Indeed, from \eqref{J-nilpotent nel th per le (n-2)}, we have that 
    \begin{equation*}
         d \varphi^{2} = B_{211} \varphi^{1 \overline{1}}, \quad d \varphi^{3} = A_{312} \varphi^{12} + B_{311} \varphi^{1 \overline{1}} + B_{312} \varphi^{1 \overline{2}} + B_{321} \varphi^{2 \overline{1}} + B_{322} \varphi^{2 \overline{2}}.
    \end{equation*}
    Thus, 
    \begin{equation*}
        \begin{split}
            &d (\overline{A_{312}} \varphi^{3 \overline{12}}) = |A_{312}|^{2} \varphi^{12 \overline{12}}, \quad d(-\overline{B_{312}} \varphi^{23 \overline{1}}) = |B_{312}|^{2} \varphi^{12 \overline{12}}, \\
            & d( -\overline{B_{321}} \varphi^{13 \overline{2}}) = |B_{321}|^{2} \varphi^{12 \overline{12}}, \quad d(\overline{B_{322}} \varphi^{13 \overline{1}}) = |B_{322}|^{2} \varphi^{12 \overline{12}}.
        \end{split}
    \end{equation*}
    Since $(M,J)$ admits a $(n-2)$-K\"ahler structure, it follows from \cite[Proposition 3.4]{HMT2023} that $A_{312} = B_{312} = B_{321} = B_{322} = 0$. Consider the left-invariant $(1,0)$-coframe $\{\psi^{1}, \dots, \psi^{n}\}$, where $\psi^{1}\doteq \varphi^{1}, \, \psi^{2}\doteq \varphi^{2} - \frac{B_{211}}{B_{311}} \varphi^{3}, \, \psi^{3} \doteq  \varphi^{3}$, and $\psi^{j} \doteq \varphi^{j}$, for $j=4, \dots, n$. Then, 
    \begin{equation*}
        d \psi^{1} = d \psi^{2} = 0.
    \end{equation*}
    By \cite[Theorem 2.3]{ST2022}, if $d \psi^{j} \neq 0$ for $j= 3, \dots, n$, then there cannot exist a $(n-2)$-K\"ahler structure. Hence, at least one $d\psi^{j} = 0$, for $j= 3, \dots, n$. Up to re-ording, $\{\psi^{1}, \dots, \psi^{n}\}$ is a left-invariant $(1,0)$-coframe that satisfies \eqref{Correct coframe}.

    \medskip
    
    Let us now consider the $d$-closed form 
    \begin{equation*}
        \Omega_{1} \doteq \sigma_{1} \big( \Omega \wedge \psi^{1 \overline{1}} + \Omega \wedge \psi^{2 \overline{2}} \big),
    \end{equation*} 
    we want to prove that $\Omega_{1}$ is transverse. 

    Let us consider $\eta \in \mathcal{A}^{1,0}(M)$ simple, defined by  
    \begin{equation*}
        \eta = \sum_{j=1}^{n} a_{j} \psi^{j},
    \end{equation*}
    where $a_{j} \in \C$, $\forall j=1, \dots, n$. We can distinguish two cases 
    \begin{enumerate}
        \item $a_{1} = 0$;
        \item $a_{1} \neq 0$.
    \end{enumerate}
    
    In the first case we have that 
    \begin{equation*}
        \sigma_{1}^{2} \, \Omega \wedge \psi^{1 \overline{1}} \wedge \eta \wedge \overline{\eta} = \sigma_{2} \Big( \Omega \wedge \eta_{1} \wedge \eta_{1}  \Big)= c_{1} \, \text{Vol},
    \end{equation*}
    where $\eta_{1} \doteq \psi^{1} \wedge \eta$ and $c_{1} > 0$ because $\Omega$ is transverse and $\eta_{1} \neq 0$, unless $a_{j}=0$, for all $j=2, \dots n$. While 
    \begin{equation*}
         \sigma_{1}^{2} \, \Omega \wedge \psi^{2 \overline{2}} \wedge \eta \wedge \overline{\eta} = \sigma_{2} \Big( \Omega \wedge \eta_{2} \wedge \eta_{2}  \Big)= c_{2} \text{Vol},
    \end{equation*}
    where $\eta_{2} \doteq \psi^{2}\wedge \eta$ and $c_{2} \geq 0$ because $\Omega$ is transverse and $\eta_{2}$ could be zero. Thus, $\sigma_{1} \, \Omega_{1} \wedge \eta \wedge \overline{\eta} > 0$.    

    While, in the second case we have that 
    \begin{equation*}
        \sigma_{1}^{2} \, \Omega \wedge \psi^{1 \overline{1}} \wedge \eta \wedge \overline{\eta} = \sigma_{2} \Omega \wedge \eta_{3} \wedge \overline{\eta_{3}} = c_{3} \, \text{Vol},
    \end{equation*}
    where $\eta_{3} \doteq \psi^{1} \wedge \eta$ and $c_{3} \geq 0$ because $\Omega$ is transverse and $\eta_{3}$ could be zero. While,
    \begin{equation*}
        \sigma_{1}^{2} \, \Omega \wedge \psi^{2 \overline{2}} \wedge \eta \wedge \overline{\eta} = \sigma_{2} \, \Omega \wedge \eta_{4} \wedge \overline{\eta_{4}}= c_{4} \text{Vol},
    \end{equation*}
    where $\eta_{4} \doteq \psi^{2} \wedge \eta$ and $c_{4}>0$ because $\Omega$ is transverse and $\eta_{4} \neq 0$. Hence, $\sigma_{1} \, \Omega_{1} \wedge \eta \wedge \overline{\eta} > 0$.
    
    Thus, $\Omega_{1}$ is a $(n-1)$-K\"ahler structure, hence $(M,J)$ is balanced. 
\end{proof}

In the following, we construct a new $3$-K\"ahler structure on the holomorphically parallelizable nilmanifold $\eta \beta_{5}$. In order to construct it, we recall the construction of transverse $(2,2)$-form on $\C^{4}$ (see \cite{BP2013, FM2025, LT2025}). 

Let us consider $\C^{4}$ endowed with the following basis of $(1,0)$-forms $\{\omega^{1}, \dots, \omega^{4}\}$. A basis for $\Lambda^{2,0}(\C^{4})$ is given by 
\begin{equation*}
    \begin{split}
        & \Omega^{1} \doteq \omega^{12}, \quad \Omega^{2} \doteq \omega^{13}, \quad \Omega^{3} \doteq \omega^{14}, \\
        &\Omega^{4} \doteq \omega^{23}, \quad \Omega^{5} \doteq - \omega^{24}, \quad \Omega^{6} \doteq \omega^{34},
    \end{split}
\end{equation*}
while, the volume form we consider, is given by $\text{Vol} \doteq \omega^{1234 \overline{1234}}$. 
\begin{remark}
    We can associate to a real $(2,2)$-form $\sigma \doteq \sum_{j,k} a_{j,k} \Omega^{j} \wedge \Omega^{\overline{k}}$ a Hermitian matrix $A_{\sigma}$ given by $A_{\sigma} = (a_{j,k})_{j,k=1,\dots,6 }$. We recall that, given $\sigma \in \Lambda^{2,2}_{\R}(\C^{4})$, we have that $\sigma$ is transverse if and only if $\overline{z}A_{\sigma}z^{t} > 0$, $\forall z\doteq(z_{1}, \dots , z_{6}) \in \C^{6}$, $z \neq 0$ such that $z \in \mathfrak{C}\doteq \{z \in \C^{6} \, | \, z_{1}z_{6} + z_{2}z_{5} + z_{3}z_{4} = 0\}$ (see \cite{BP2013, FM2025, LT2025}).
    
    Furthermore, it can be easily seen that 
    \begin{equation*}
        \Omega^{j} \wedge \Omega^{k} = \begin{cases}
            \omega^{1234}, \quad \text{if} \, \,  k = 7-j, \\
            0, \quad \text{otherwise}.
        \end{cases}
    \end{equation*}
\end{remark}
\begin{lemma}\label{a new transverse form in C^{4}}
    The $(2,2)$-form 
    \begin{equation}
        \begin{split}
            \Omega_{a,b} = & \, 4 \Omega_{1} \wedge \overline{\Omega_{1}} + \Omega_{2} \wedge \overline{\Omega_{2}} + 4\Omega_{3} \wedge \overline{\Omega_{3}} + \Omega_{4} \wedge \overline{\Omega_{4}} + \Omega_{5} \wedge \overline{\Omega_{5}}   \\
            & + \Omega_{6} \wedge \overline{\Omega_{6}} + a \,  \Omega_{2} \wedge \overline{\Omega_{5}} + \overline{a} \, \Omega_{5} \wedge \overline{\Omega_{2}} + b \, \Omega_{3} \wedge \overline{\Omega_{4}} + \overline{b} \, \Omega_{4} \wedge \overline{\Omega_{3}},
        \end{split}
    \end{equation}
    for $a,b \in \C$, $|a|, |b|\leq 1$, is transverse.
\end{lemma}
\begin{proof}
    Let $A$ be the matrix associated to $\Omega_{a,b}$ and let $z \in \mathfrak{C}\doteq \{z \in \C^{6} \, | \, z_{1}z_{6} + z_{2}z_{5} + z_{3}z_{4} = 0\}$. We have that 
    \begin{equation*}
        \begin{split}
            \overline{z}A z & = 4 |z_{1}|^{2} + |z_{2}|^{2} + 4 |z_{3}|^{2} + |z_{4}|^{2} + |z_{5}|^{2} + |z_{6}|^{2} + 2 \text{Re} (a \overline{z_{2}} z_{5}) + 2 \text{Re}(b \overline{z_{3}}z_{4}) \\
             & \geq 4 |z_{1}z_{6}| + 2|z_{2}z_{5}| + 4|z_{3}z_{4}| + 2 \text{Re} (a \overline{z_{2}} z_{5}) + 2 \text{Re}(b \overline{z_{3}}z_{4}) \\
             & \geq 4 |z_{1}z_{6}| + 2|z_{2}z_{5}| + 4|z_{3}z_{4}| - 2|a  z_{2} z_{5}| - 2 |b  z_{3} z_{4}| \\
             & \geq 4 |z_{1}z_{6}| + 2|z_{2}z_{5}| + 4|z_{3}z_{4}| - 2|z_{2} z_{5}| - 2 |z_{3}z_{4}| \\
             & \geq 2\big(2 |z_{1}z_{6}| + |z_{3}z_{4}| \big) \geq 0.
        \end{split}
    \end{equation*}
    We notice that, the inequalities are equalities if and only if 
    \begin{equation*}
        \begin{split}
            & 2|z_{1}| = |z_{6}|, \quad |z_{2}| = |z_{5}|, \quad 2|z_{3}| = |z_{4}|, \\
            & a\overline{z_{2}}z_{3} \in \R_{\leq 0}, \quad a\overline{z_{3}}z_{4} \in \R_{\leq 0}, \\
            & |z_{1}z_{6}| = |z_{3}z_{4}| = 0.
        \end{split}
    \end{equation*}
    Hence, for the equalities we must have $z_{1}=z_{6}=z_{3}=z_{4}=0$. Furthermore, since $z \in \mathfrak{C}$, also $z_{2}=z_{5}=0$. Thus, $\Omega_{a,b}$ is transverse.
\end{proof}
Now, we are ready to construct a new $3$-K\"ahler structure on $\eta \beta_{5}$.
\begin{example}\label{etabeta_{5}}
    We recall that, according to the classification of I. Nakamura (see \cite{N1975}), $\eta \beta_{5} \doteq \Gamma \backslash G$, where $G$ can be identified with $(\C^{5}, \ast)$ and the multiplication is defined by
    \begin{equation*}
        (y_{1}, \dots, y_{5}) \ast (z_{1},\dots, z_{5}) \doteq (y_{1} + z_{1}, \dots, y_{5} + z_{5} + y_{1}z_{3} + y_{2}z_{4}),
    \end{equation*} 
    and $\Gamma$ is the lattice composed by Gaussian integers. A
    basis of left-invariant $(1,0)$-forms is given by 
    \begin{equation*}
        \varphi^{1} \doteq dz^{1}, \quad \varphi^{2} \doteq dz^{2}, \quad \varphi^{3} \doteq dz^{3}, \quad \varphi^{4} \doteq dz^{4}, \quad  \varphi^{5} = dz^{5}- z^{1} dz^{3} - z^{2}dz^{4},
    \end{equation*}
    where $(z^{1}, \dots, z^{5})$ are the holomorphic coordinates, meanwhile the dual frame of left-invariant vector fields of type $(1,0)$ is given by 
    \begin{equation*}
        \theta_{1} \doteq \frac{\del}{\del z^{1}}, \quad \theta_{2} \doteq \frac{\del}{\del z^{2}}, \quad \theta_{3} \doteq \frac{\del}{\del z^{3}} + z^{1}\frac{\del}{\del z^{5}} , \quad \theta_{4} \doteq \frac{\del}{\del z^{4}} +  z^{2}\frac{\del}{\del z^{5}}, \quad \theta_{5} \doteq \frac{\del}{\del z^{5}}.
    \end{equation*}
    Furthermore,
    \begin{equation*}
        d \varphi^{1} = \dots = d \varphi^{4} = 0, \quad d \varphi^{5} = - \varphi^{13} - \varphi^{24}.
    \end{equation*}
    
    We recall that $\eta \beta_{5}$ is balanced and it admits a $3$-K\"ahler structure (see \cite{AB1991}). The $3$-K\"ahler structure constructed in \cite{AB1991} is given by 
    \begin{equation}\label{old 3-kahler form}
        \Omega \doteq \sigma_{3} \big( \sum_{i<j<k} \varphi^{ijk \overline {ijk}} - \varphi^{135 \overline{245}} - \varphi^{245 \overline{135}} \big).
    \end{equation}

    Now, let us consider 
    \begin{equation}\label{new 3-kahler form}
        \begin{split}
            \Omega_{\ast} \doteq & \sigma_{3} \big( \sum_{i<j<k<5} \varphi^{ijk \overline{ijk}} + 4\varphi^{125 \overline{125}} + \varphi^{135 \overline{135}} + 4 \varphi^{145 \overline{145}} + \varphi^{235 \overline{235}} \\
            & + \varphi^{245 \overline{245}} + \varphi^{345 \overline{345}} - \varphi^{135 \overline{245}} - \varphi^{245 \overline{135}} - \varphi^{235 \overline{145}} - \varphi^{145 \overline{235}} \big).
        \end{split}
    \end{equation}
    An easy computation shows that $d\Omega_{\ast}=0$. 

    To show that $\Omega_{\ast}$ is transverse, we make use of \cite[Theorem 3.5]{FM2025}. Let us consider $\{v_{1}, \dots, v_{5}\}$, where $v_{j} \doteq \theta_{j}|_{e}$, for $j=1, \dots , 5$ and $e$ is the identity of $(G,\ast)$. Decompose $\mathfrak{g_{\C}}$, which is the complexified Lie algebra of $\eta\beta_{5}$ in the following way 
    \begin{equation*}
        \mathfrak{g}_{\C} = \mathfrak{h}_{\C} \oplus \langle v_{5}, \overline{v_{5}} \rangle,
    \end{equation*}
    where $\mathfrak{h}_{\C} \doteq \mathfrak{h}_{1,0} \oplus \overline{\mathfrak{h}_{1,0}}$ and $\mathfrak{h}_{1,0} \doteq \langle v_{1}, \dots, v_{4}\rangle$. 

    From the definition of $\Omega_{\ast}$ we have that 
    \begin{equation*}
        \Omega_{\ast} =   \Omega_{\ast}|_{\mathfrak{h}_{\C}} + \sigma_{3} \, \varphi^{5 \overline{5}} \wedge F,
    \end{equation*}
    where 
    \begin{equation*}
        \begin{split}
            & \Omega_{\ast}|_{\mathfrak{h}_{\C}} \doteq \sigma_{3} \sum_{i<j<k<5} \varphi^{ijk \overline{ijk}}, \\
            & F \doteq 4\varphi^{12 \overline{12}} + \varphi^{13 \overline{13}} + 4 \varphi^{14 \overline{14}} + \varphi^{23 \overline{23}} + \varphi^{24 \overline{24}} + \varphi^{34 \overline{34}} - \varphi^{13 \overline{24}} - \varphi^{24 \overline{13}} - \varphi^{23 \overline{14}} - \varphi^{14 \overline{23}}.
        \end{split}
    \end{equation*}
    We can easily see that $\Omega_{\ast}|_{\mathfrak{h}_{\C}}$ is a tranverse $(3,3)$-form on $\mathfrak{h}_{\C}$, indeed it is the third power of the fundamental form of the diagonal Hermitian metric on $\mathfrak{h}_{\C}$. Whereas, by Lemma \ref{a new transverse form in C^{4}}, we have that $F$ is a transverse $(2,2)$-form on $\mathfrak{h}_{\C}$. Thus, by \cite[Theorem 3.5]{FM2025}, $\Omega_{\ast}$ is transverse. 

    As required, $\Omega_{\ast}$ is a new example of $3$-K\"ahler structure on $\eta\beta_{5}$.
\end{example}

In the following example, we construct a family of $3$-K\"ahler nilmanifolds that is not holomorphically parallelizable. 

\begin{example}\label{3-kahler nilmanifolds}
    Let us consider the set of $(1,0)$-forms $\{ \sigma^{1}, \dots, \sigma^{5}\}$, such that 
    \begin{equation}\label{3-Kahler nilmanifolds examples}
        \begin{cases}
            d \sigma^{j}=0, \quad \forall j= 1, \dots, 4, \\
            d \sigma^{5} = a \, \sigma^{13} + b \, \sigma^{1 \overline{2}} + c \, \sigma^{1 \overline{4}} + a \, \sigma^{24} + d \, \sigma^{2 \overline{1}} \\
            \quad \quad \; + c \, \sigma^{2 \overline{3}} + e \, \sigma^{3 \overline{2}} - d \, \sigma^{3 \overline{4}} + e \, \sigma^{4 \overline{1}} - b \, \sigma^{4 \overline{3}},
        \end{cases}
    \end{equation}
    where $a,b,c,d,e \in \Q[i]$. Then, the complex forms $\{\sigma^{1}, \dots, \sigma^{5}\}$ span the dual of a Lie algebra $\mathfrak{g}$, which is $2$-step nilpotent and depends on the choice of the coefficients. Furthermore, the almost complex structure $J$ defined by \eqref{3-Kahler nilmanifolds examples} is integrable and nilpotent. Let us denote by $G$ the simply connected, connected Lie group whose Lie algebra is denoted by $\mathfrak{g}$. Due to Malcev's Theorem \cite{M1951} there exists a lattice $\Gamma$ such that $(M \doteq \Gamma \backslash G,J)$ is a nilmanifold, where $J$ is the complex structure induced by the complex structure on the Lie algebra $\mathfrak{g}$. 

    We recall that this family of nilmanifolds is a particular case of the family of nilmanifolds constructed in \cite{ST2023}.

    Let us consider the real $(3,3)$-form 
    \begin{equation*}
        \Omega \doteq \sigma_{3} \Big(\sum_{i <j < k} \sigma^{ijk \overline{ijk}} - \sigma^{135 \overline{245}} - \sigma^{245 \overline{135}} \Big).
    \end{equation*}
    As shown in \cite{AB1991, LT2025}, it is a transverse form, moreover it is also $d$-closed. Indeed, 
    \begin{equation*}
        \begin{split}
            &\delbar \sigma^{ijk \overline{ijk}}= 0, \quad i,j,k \neq 5, \\
            & \delbar \sigma^{125 \overline{125}} = - d \sigma^{123 \overline{1245}} - b \sigma^{124 \overline{1235}}, \quad \delbar \sigma^{135 \overline{135}} =  \overline{a} \sigma^{135 \overline{1234}}, \\
            & \delbar \sigma^{145 \overline{145}} = c \sigma^{124 \overline{1345}} + e \sigma^{134 \overline{1245}}, \quad \delbar \sigma^{235 \overline{235}} = c \sigma^{123 \overline{2345}} + e \sigma^{234 \overline{1235}}, \\
            & \delbar \sigma^{345 \overline{345}} = b \sigma^{134 \overline{2345}} + d \sigma^{234 \overline{1345}}, \quad  \delbar \sigma^{245 \overline{245}} = \overline{a} \sigma^{245 \overline{1234}},
        \end{split}
    \end{equation*}
    and
    \begin{equation*}
        \begin{split}
            & \delbar \sigma^{135 \overline{245}} = - d \sigma^{123 \overline{1245}} + c \sigma^{123 \overline{2345}} + e \sigma^{134 \overline{1245}} + b \sigma^{134 \overline{2345}} + \overline{a} \sigma^{135 \overline{1234}}, \\
            & \delbar \sigma^{245 \overline{135}} = - b \sigma^{124 \overline{1235}} + c \sigma^{124 \overline{1345}} + e \sigma^{234 \overline{1235}} + d \sigma^{234 \overline{1345}} + \overline{a} \sigma^{245 \overline{1234}},
        \end{split}
    \end{equation*}
    thus, $\delbar \Omega=0$. Hence, $\Omega$ is a $3$-K\"ahler structure on $(M,J)$.
\end{example}

We conclude this section by providing a family of $3$-symplectic nilmanifolds of complex dimension $5$.
\begin{example}\label{3-symplectic nilmanifolds}
    Let $\{\sigma^{1}, \dots, \sigma^{5}\}$ be the set of complex $(1,0)$-forms such that 
    \begin{equation}\label{Structure equations for 10-dimensional Nilmanifolds with Astheno and 2-pluriclosed}
        \begin{cases}
            d \sigma^{j}=0, \quad \forall j= 1, \dots, 4, \\
            d \sigma^{5} = a \, \sigma^{12} + c \, \sigma^{34}+ b \, (\sigma^{1 \overline{1}} +  \sigma^{2 \overline{2}} +  \sigma^{3 \overline{3}} +  \sigma^{4 \overline{4}}) ,
        \end{cases}
    \end{equation}
    where $a,b,c, \in \Q[i]$. Then, the complex forms $\{\sigma^{1}, \dots, \sigma^{5}\}$ span the dual of a $2$-step nilpotent Lie algebra $\mathfrak{g}$ which depends on the choice of the coefficients. Furthermore, the almost complex structure $J$ defined by \eqref{Structure equations for 10-dimensional Nilmanifolds with Astheno and 2-pluriclosed} is integrable and nilpotent. Let us denote by $G$ the simply connected, connected Lie group whose Lie algebra is $\mathfrak{g}$. We recall that Malcev's Theorem \cite{M1951} implies that there exists a lattice $\Gamma$ such that $(M \doteq \Gamma \backslash G,J)$ is a nilmanifold, where $J$ is the complex structure induced by the complex structure on the Lie algebra $\mathfrak{g}$. We mention that also this family of nilmanifolds is a particular case of the family constructed in \cite{ST2023}.
    
    From now on, we assume that $a,b,c \neq 0$. Thus, by a straightforward computation using \cite[Proposition 3.4]{HMT2023}, we conclude that there are no $p$-K\"ahler structures for $p=2,3,4$. Moreover, since $M = \Gamma \backslash G$ and $G$ is a non-abelian Lie group, there are no K\"ahler metrics on $M$. 
    
    Let us consider the diagonal left-invariant metric $g$ and denote by $\omega$ the fundamental form of $g$. Define 
    \begin{equation*}
        \Omega \doteq \, \omega^{3} + \beta_{1} + \overline{\beta_{1}} + \beta_{2} + \overline{\beta_{2}}, 
    \end{equation*}
    where $\beta_{1} \doteq L \sigma^{1245 \overline{45}} + M \sigma^{1235 \overline{35}} + N \sigma^{1345 \overline{15}} + O \sigma^{2345 \overline{25}}$, for $L,M,N,O \in \C$, $\beta_{2} \doteq P \sigma^{12345 \overline{5}}$, for $P \in \C$. Suppose that $L, M, N, O, P \neq 0$. It can be easily seen that the real $(3,3)$-form $\Omega$ is a $3$-symplectic structure if and only if 
    \begin{equation*}
        \del \omega^{3} = \delbar \beta_{1}, \quad \del \beta_{1} = \delbar \beta_{2}.
    \end{equation*}
    Since
    \begin{equation*}
        \begin{split}
            \omega^{3} = & \,\frac{3}{4} i(\sigma^{123 \overline{123}} + \sigma^{124 \overline{124}} + \sigma^{125 \overline{125}} + \sigma^{134 \overline{134}} + \sigma^{135 \overline{135}} \\
            & + \sigma^{145 \overline{145}} + \sigma^{234 \overline{234}} + \sigma^{235 \overline{235}} + \sigma^{245 \overline{245}} + \sigma^{345 \overline{345}}),
        \end{split}
    \end{equation*}
    then 
    \begin{equation*}
        \begin{split}
            \del \omega^{3} = \, \frac{3}{4} i \Big(- 3 \overline{b} \, \big( \sigma^{1235 \overline{123}} + \sigma^{1245 \overline{124}} + \sigma^{1345 \overline{134}} + \sigma^{2345 \overline{234}} \big) + a \, \sigma^{1234 \overline{345}} + c \, \sigma^{1234 \overline{125}}\Big).
        \end{split}
    \end{equation*}
    Furthermore, 
    \begin{equation*}
        \begin{split}
            \delbar \beta_{1} = & \, b \Big(L + M \Big) \sigma^{1234 \overline{345}} + b \Big(N + O  \Big) \sigma^{1234 \overline{125}} \\
            & - \overline{a} L \sigma^{1245 \overline{124}} - \overline{a} M \sigma^{1235 \overline{123}} - \overline{c} N \sigma^{1345 \overline{134}} - \overline{c} O \sigma^{2345 \overline{234}},
        \end{split}
    \end{equation*}
    so, in order to have $\del \omega^{3} = \delbar \beta_{1}$, we have that 
    \begin{equation*}
        \begin{cases}
            b (L+M) = \frac{3}{4}i \,a, \\
            b (N+O) = \frac{3}{4}i \,c, \\
            \overline{a} \, L = \frac{9}{4}i \, \overline{b}, \\
            \overline{a} \, M = \frac{9}{4}i \, \overline{b}, \\
            \overline{c} \, N = \frac{9}{4}i \, \overline{b}, \\
            \overline{c} \, O = \frac{9}{4}i \, \overline{b}.
        \end{cases}
    \end{equation*}
    The last $4$ equations imply that $L = M$ and $N = O$. From now on, we assume that $L = M$ and $N = O$, thus the system becomes
    \begin{equation*}
        \begin{cases}
            2 b L  = \frac{3}{4}i \, a, \\
            2 b N = \frac{3}{4}i \, c, \\
            \overline{a} L = \frac{9}{4}i \, \overline{b}, \\
            \overline{c} N = \frac{9}{4}i \, \overline{b}.
        \end{cases}
    \end{equation*}
    By simple calculations, we get that the solution of the system is given by
    \begin{equation}\label{condition for L and a1}
        L = i \frac{3a}{8b}, \quad N = i\frac{3c}{8 b}, \quad |a|^{2} = |c|^{2} = 6 |b|^{2}.
    \end{equation}
    If the parameters $L, N, a, b, c$ satisfy \eqref{condition for L and a1}, we get that 
    \begin{equation*}
        \del \omega^{3} = \delbar \beta_{1}.
    \end{equation*}
    Since 
    \begin{equation*}
        \delbar \beta_{2} = - \overline{a} P \sigma^{12345 \overline{12}} - \overline{c} P \sigma^{12345 \overline{34}}, \quad \del \beta_{1} = 2 \overline{b} L \sigma^{12345 \overline{34}} + 2 \overline{b} N  \sigma^{12345 \overline{12}},
    \end{equation*}
    we have that $\del \beta_{1} = \delbar \beta_{2}$ if and only if 
    \begin{equation}\label{condition for P}
        \begin{cases}
            P = - \frac{2 \overline{b}L}{\overline{c}}, \\
            P = - \frac{2\overline{b}N}{\overline{a}}.
        \end{cases}
    \end{equation}
    Note that, if \eqref{condition for L and a1} holds, then $\frac{L}{\overline{c}} = \frac{N}{\overline{a}}$.
    
    Thus, if 
    \begin{equation*}
        \beta_{1} \doteq \frac{3}{8b}i \Big( a \big( \sigma^{1245 \overline{45}} + \sigma^{1235 \overline{35}} \big) + c \big( \sigma^{1345 \overline{15}} +  \sigma^{2345 \overline{25}} \big) \Big), \quad \beta_{2} \doteq - \frac{3\overline{b}a}{4b\overline{c}}i \, \sigma^{12345 \overline{5}},
    \end{equation*}
    and $|a|^{2} = |c|^{2} = 6 |b|^{2}$, then
    \begin{equation*}
        \Omega = \omega^{3} + \beta_{1} + \overline{\beta_{1}} + \beta_{2} + \overline{\beta_{2}}
    \end{equation*}
    is a $3$-symplectic structure on $(M,J)$. Moreover, by \cite[Theorem 4.1]{ST2023}, $\omega$ is an astheno-K\"ahler metric and by \cite[Theorem 6.9]{LT2025}, $(M,J)$ admits a $4$-symplectic structure. 
\end{example}

\section{Deformations of p-K\"ahler structures}\label{Section 5}
In this section we provide obstruction results for the existence of smooth curves of $p$-K\"ahler structures along a differentiable family of compact complex manifolds. 

Let $\{M_{t}\}_{t \in I}$, where $I \doteq (-\epsilon, \epsilon)$, for $\epsilon >0$, be a curve of deformations of $(M,J)$, which is a compact complex manifold of complex dimension $n$. Suppose that the curve of deformations $\{M_{t}\}_{t \in I}$ is parametrized by the $(0,1)$-vector form $\varphi(t) \in \mathcal{A}^{0,1}(T^{1,0}M)$. Moreover, suppose that $(M,J)$ admits a $p$-K\"ahler structure denoted by $\Omega$, for $1<p<n-1$. 

Let $\{\Omega_{t}\}_{t \in I}$ be a smooth family of real transverse $(p,p)$-forms along $\{M_{t}\}$ such that $\Omega_{0} \doteq \Omega$. Since the operator $e^{\iota_{\varphi(t)}|\iota_{\overline{\varphi(t)}}}$ is a linear isomorphism, then $\{\Omega_{t}\}_{t \in I}$ can be written as 
\begin{equation}\label{smooth family of (p,p)-forms}
    \Omega_{t} = e^{i_{\varphi}|i_{\overline{\varphi}}} \big(\Omega(t) \big),
\end{equation}
where, locally $\Omega(t) = \Omega_{I_{p}J_{p}}(t) dz^{I_{p}}\wedge dz^{\overline{J_{p}}} \in \mathcal{A}^{p,p}(M)$.

The following theorem is proved with the same techniques adopted in \cite{PS2021, S2022, S2023}. For the sake of completeness, we outline the proof in our case.  
\begin{theorem}\label{Necessary condition for the existence of curve of 3-K}
    Let $(M,J,\Omega)$ be a compact $p$-K\"ahler manifold of complex dimension $n$, for $1 < p < n-1$, and let $\{M_{t}\}_{t \in I}$ be a differentiable family of compact complex manifolds, where $I \doteq (-\epsilon, \epsilon), \epsilon >0$ and $M_{0} \doteq M$. Suppose that $\{M_{t}\}_{t \in I}$ is parametrized by the $(0,1)$-vector form $\varphi(t) \in \mathcal{A}^{0,1}(T^{1,0} M)$. 

    Let $\{\Omega_{t}\}_{t \in I}$ be a smooth family of real transverse $(p,p)$-forms along $\{M_{t}\}_{t \in I}$ as in \eqref{smooth family of (p,p)-forms}. If every $\Omega_{t}$ is a $p$-K\"ahler structure, then
    \begin{equation}
        \del \circ \iota_{\varphi^{'}(0)}(\Omega) = - \delbar \, \Omega^{'}(0).
    \end{equation}
\end{theorem}

\begin{proof}
    If every $\Omega_{t}$ is a $p$-K\"ahler structure, then $\delbar_{t} \Omega_{t} = 0$, i.e., 
    \begin{equation*}
        \delbar_{t} \Big( e^{\iota_{\varphi}|\iota_{\overline{\varphi}}}  \big(\Omega(t) \big) \Big) = 0.
    \end{equation*}
    By means of \eqref{delbar operator on M_t}, we have that 
    \begin{equation*}
        \delbar_{t} \Big( e^{\iota_{\varphi}|\iota_{\overline{\varphi}}}  \big(\Omega(t) \big) \Big) = e^{\iota_{\varphi}|\iota_{\overline{\varphi}}} \Big(\big(I - \overline{\varphi} \varphi \big)^{-1} \Finv \, \big([\del,\iota_{\varphi}] + \delbar \big) \big(I - \overline{\varphi} \varphi \big) \Finv \, \Omega(t) \Big),
    \end{equation*}
    which can be rewritten as 
    \begin{equation*}
        \delbar_{t} \Big( e^{\iota_{\varphi}|\iota_{\overline{\varphi}}}  \big(\Omega(t) \big) \Big) = (I + \varphi + \overline{\varphi}) \, \Finv \, \Big(\big(I - \overline{\varphi} \varphi \big)^{-1} \Finv \,\big([\del,\iota_{\varphi}] + \delbar \big) (I - \overline{\varphi} \varphi \big) \Finv \, \Omega(t) \Big).
    \end{equation*}
    By using the Taylor series expansion centered in $t=0$ for $\varphi(t)$, we get $\varphi(t)= t \varphi^{'}(0) + o(t)$. 
    
    Therefore, 
    \begin{equation*}
        \big(I - \varphi \overline{\varphi} \big) = \big(I - \overline{\varphi} \varphi \big) = \big(I - \varphi \overline{\varphi} \big)^{-1} = \big(I - \overline{\varphi} \varphi \big)^{-1} = I + o(t), 
    \end{equation*}
    and by using the Taylor series expansion centered in $t=0$ for $\Omega(t)$, we get that 
    \begin{equation*}
        \begin{split}
            \delbar_{t} \Omega_{t} & = (I + t \varphi^{'}(0) + t \overline{\varphi^{'}(0)}) \, \Finv \, \Big( \big([\del,t \,\varphi^{'}(0) \lrcorner] + \delbar \big) \big(\Omega(0) + t \, \Omega^{'}(0) \big) \Big) + o(t) \\
            & = (I + t \varphi^{'}(0) + t \overline{\varphi^{'}(0)}) \, \Finv \, \Big( t \del \big(\varphi^{'}(0) \lrcorner \, \Omega(0) \big) + t  \delbar \, \Omega^{'}(0)  \Big) + o(t) \\
            & = t \del \big(\varphi^{'}(0) \lrcorner \, \Omega(0) \big) + t  \delbar \, \Omega^{'}(0) + o(t).
        \end{split}     
    \end{equation*}
    Since $\delbar_{t} \Omega_{t} = 0$, $\forall t \in I$, differentiating, we get 
    \begin{equation*}
        \begin{split}
            \frac{\del}{\del  t} (\delbar_{t} \Omega_{t}) & = \frac{\del}{\del  t} \Big( t \del \big(\varphi^{'}(0) \lrcorner \, \Omega(0) \big) + t  \delbar \, \Omega^{'}(0) + o(t) \Big) \\
            & = \del \big(\varphi^{'}(0) \lrcorner \, \Omega(0) \big) + \delbar \, \Omega^{'}(0) = 0.
        \end{split}
    \end{equation*}
    As required.
\end{proof}
A direct consequence of Theorem \ref{Necessary condition for the existence of curve of 3-K} is the following.
\begin{corollary}\label{Obstruction in dolbeault cohomology}
    Let $(M,J,\Omega)$ be a compact $p$-K\"ahler manifold of complex dimension $n$, for $1 < p < n-1$, and let $\{M_{t}\}_{t \in I}$ be a differentiable family of compact complex manifolds, where $I \doteq (-\epsilon, \epsilon), \epsilon >0$ and $M_{0} \doteq M$. Suppose that $\{M_{t}\}_{t \in I}$ is parametrized by the $(0,1)$-vector form $\varphi(t) \in \mathcal{A}^{0,1}(T^{1,0} M)$. If there exists a smooth family of $p$-K\"ahler structures $\{\Omega_{t}\}_{t \in I}$ such that $\Omega_{0} \doteq \Omega$, then 
    \begin{equation*}
        \Big[ \del \circ \iota_{\varphi^{'}(0)} (\Omega) \Big]_{H^{p,p+1}_{\delbar}(M)} = 0.
    \end{equation*}
\end{corollary}
In the following, we study some examples of curves of deformations of $\eta \beta_{5}$. 
\begin{example}\label{examples of deformations}
    Let us consider Example \ref{etabeta_{5}}. We recall that in \cite{AB1990}, the authors prove that, in general, the property of being $p$-K\"ahler is not stable under small deformations. The deformation studied by the authors of \cite{AB1990} is given by
    \begin{equation*}
        \Psi(t) \doteq t \, \varphi^{\overline{1}} \otimes \theta_{1}.
    \end{equation*}
    Here, we consider the small deformations of $\eta \beta_{5}$ parametrized by the $(0,1)$-vector form $\Psi(\mathbf{t}) \in \mathcal{A}^{0,1}\big(T^{1,0}(\eta \beta_{5})\big)$ defined by
    \begin{equation*}
        \Psi(\mathbf{t})= t_{1} \, \varphi^{\overline{1}} \otimes \theta_{1} + t_{2} \, \varphi^{\overline{2}} \otimes \theta_{2},
    \end{equation*}
    where $\mathbf{t} \doteq (t_{1},t_{2}) \in B(0,\delta)\subset \C^{2}$, $\delta > 0$. 
    
    By solving 
    \begin{equation*}
        \begin{cases}
            \big(\delbar - \Psi(\mathbf{t}) \big) \zeta^{j}_{\mathbf{t}} = 0, \\
            \zeta^{j}_{0} = z^{j}, \quad j = 1, \dots , 5,
        \end{cases}
    \end{equation*}
    where $(z^{1}, \dots, z^{5})$ are the holomorphic coordinates on $\eta \beta_{5}$, we get that 
    \begin{equation*}
        \begin{cases}
            \zeta^{1}_{\mathbf{t}} = z^{1} + t_{1} \overline{z^{1}} \\
            \zeta^{2}_{\mathbf{t}} = z^{2} + t_{2} \overline{z^{2}} \\
            \zeta^{j}_{\mathbf{t}} = z^{j}, \quad j = 3,4,5.
        \end{cases}
    \end{equation*}
    Therefore, a basis of left-invariant $(1,0)$-form is given by 
    \begin{equation*}
        \begin{cases}
            \varphi^{j}_{\mathbf{t}} = d \zeta^{j}_{\mathbf{t}}, \quad j =1, \dots, 4, \\
            \varphi^{5}_{\mathbf{t}} = d \zeta^{5}_{\mathbf{t}} + \frac{1}{1 - |t_{1}|^{2}}(t_{1}\zeta^{\overline{1}}_{\mathbf{t}} - \zeta^{1}_{\mathbf{t}}) d \zeta^{3}_{\mathbf{t}} + \frac{1}{1 - |t_{2}|^{2}} (t_{2} \zeta^{\overline{2}}_{\mathbf{t}} - \zeta^{2}_{\mathbf{t}}) d \zeta^{4}_{\mathbf{t}}.
        \end{cases}
    \end{equation*}
    The complex structure equations for the coframe $\{\varphi^{1}_{\mathbf{t}}, \dots, \varphi^{5}_{\mathbf{t}}\}$ are given by 
    \begin{equation*}
        \begin{cases}
            d \varphi^{j}_{\mathbf{t}} = 0, \quad j=1, \dots 4, \\
            d \varphi^{5}_{\mathbf{t}} = - \frac{t_{1}}{1 - |t_{1}|^{2}} \varphi^{3}_{\mathbf{t}} \wedge \varphi^{\overline{1}}_{\mathbf{t}} - \frac{1}{1 - |t_{1}|^{2}} \varphi^{1}_{\mathbf{t}} \wedge \varphi^{3}_{\mathbf{t}} - \frac{t_{2}}{1 - |t_{2}|^{2}} \varphi^{4}_{\mathbf{t}} \wedge \varphi^{\overline{2}}_{\mathbf{t}} - \frac{1}{1 - |t_{2}|^{2}} \varphi^{2}_{\mathbf{t}} \wedge \varphi^{4}_{\mathbf{t}}. 
        \end{cases}
    \end{equation*}
    Note that the almost complex structure defined on each element of the deformation is integrable. 
    
    Let us consider the smooth curve of deformations given by 
    \begin{equation*}
        t \mapsto \Psi(t) \doteq t (a_{1} \varphi^{\overline{1}} \otimes \theta_{1} + a_{2} \varphi^{\overline{2}} \otimes\theta_{2}), \quad t \in I =(-\epsilon, \epsilon), \; \epsilon > 0, \; a_{1}, a_{2} \in \C.
    \end{equation*}
    We can easily see that the derivative of $\Psi(t)$ in $0$ is given by 
    \begin{equation*}
        \Psi^{'}(0) = a_{1} \varphi^{\overline{1}} \otimes \theta_{1} + a_{2} \varphi^{\overline{2}} \otimes\theta_{2}.
    \end{equation*}
    
    Suppose that $\{\Omega_{t}\}_{t \in I}$ and $\{\widetilde{\Omega}_{t}\}_{t \in I}$ are two smooth families of $3$-K\"ahler structures along $\{M_{t}\}_{t \in I}$ such that $\Omega_{0} = \Omega$ and $\widetilde{\Omega}_{t} = \Omega_{\ast}$, where $\Omega$ and $\Omega_{\ast}$ are defined as in \eqref{old 3-kahler form} and \eqref{new 3-kahler form}, respectively. 

    \medskip
    
    In the first case we have that 
    \begin{equation*}
        \del \circ \iota_{\Psi^{'}(0)}(\Omega) = 
        \sigma_{3}\big( a_{1} \varphi^{234 \overline{1245}} + a_{2} \varphi^{134 \overline{1235}} \big).
    \end{equation*}
    We can easily see that the classes $[\varphi^{234 \overline{1245}}]_{\delbar}$ and $[\varphi^{134 \overline{1235}}]_{\delbar}$ are non vanishing since $\varphi^{234 \overline{1245}}$ and $\varphi^{134 \overline{1235}}$ are harmonic with respect to the Dolbeault Laplacian. Hence, by Corollary \ref{Obstruction in dolbeault cohomology} there exists no curve of $3$-K\"ahler structures $\{\Omega_{t}\}_{t \in I}$ such that $\Omega_{0} = \Omega$ along the curve of deformations $t \mapsto \Psi(t)$.

    \medskip
    
    Meanwhile, in the second case, we have that 
    \begin{equation*}
        \begin{split}
            \del \circ \iota_{\Psi^{'}(0)}(\Omega_{\ast}) & = \sigma_{3} \del \Big(- a_{1}\varphi^{35 \overline{1245}} + a_{2}\varphi^{45 \overline{1235}} + a_{2} \varphi^{35\overline{1245}} - a_{1} \varphi^{45 \overline{1235}} \Big) \\
            & = \sigma_{3} \Big( (a_{1} - a_{2}) \varphi^{234 \overline{1245}} + (a_{2} - a_{1}) \varphi^{134 \overline{1235}} \Big).
        \end{split}
    \end{equation*}
    By using the same argument as before, there exists no curve of $3$-K\"ahler structures $\{\widetilde{\Omega}_{t}\}_{t \in I}$ such that $\Omega_{0} = \Omega_{\ast}$ along the curve of deformations $t \mapsto \Psi(t)$ unless $a_{1} = a_{2}$. If $a_{1} = a_{2}$, then we can construct a curve of $3$-K\"ahler structures. 
    
    If $a_{1} = a_{2}$, we have that the holomorphic coordinates along the deformation are given by 
    \begin{equation*}
        \begin{cases}
            \zeta^{1}_{t} = z^{1} + t a_{1} \overline{z^{1}}, \\
            \zeta^{2}_{t} = z^{2} + t a_{1} \overline{z^{2}}, \\
            \zeta^{j}_{t} = z^{j}, \quad j = 3,4,5.
        \end{cases}
    \end{equation*}
    Therefore, a basis of left-invariant $(1,0)$-form is given by 
    \begin{equation*}
        \begin{cases}
            \varphi^{j}_{t} = d \zeta^{j}_{t}, \quad j =1, \dots, 4, \\
            \varphi^{5}_{t} = d \zeta^{5}_{t} + T (t a_{1} \zeta^{\overline{1}}_{t} - \zeta^{1}_{t}) d \zeta^{3}_{t} + T ( t a_{1}\zeta^{\overline{2}}_{t} - \zeta^{2}_{t}) d \zeta^{4}_{t},
        \end{cases}
    \end{equation*}
    where $T \doteq \frac{1}{1- |ta_{1}|^{2}}$. Hence, the structure equations for the coframe $\{\varphi^{1}_{t}, \dots, \varphi^{5}_{t}\}$ are given by 
    \begin{equation*}
        \begin{cases}
            d \varphi^{j}_{t} = 0, \quad j=1, \dots 4, \\
            d \varphi^{5}_{t} = - T t a_{1} \varphi^{3}_{t} \wedge \varphi^{\overline{1}}_{t} - T \varphi^{1}_{t} \wedge \varphi^{3}_{t} - T t a_{1} \varphi^{4}_{t} \wedge \varphi^{\overline{2}}_{t} - T \varphi^{2}_{t} \wedge \varphi^{4}_{t}, 
        \end{cases}
    \end{equation*}
    Let us consider the following family of $(3,3)$-forms $\{\widetilde{\Omega}_{t}\}$, defined by
    \begin{equation*}
        \begin{split}
            \widetilde{\Omega}_{t} \doteq & \sigma_{3} \big( \sum_{i<j<k<5} \varphi_{t}^{ijk \overline{ijk}} + 4\varphi_{t}^{125 \overline{125}} + \varphi_{t}^{135 \overline{135}} + 4 \varphi_{t}^{145 \overline{145}} + \varphi_{t}^{235 \overline{235}} \\
            & + \varphi_{t}^{245 \overline{245}} + \varphi_{t}^{345 \overline{345}} - \varphi_{t}^{135 \overline{245}} - \varphi_{t}^{245 \overline{135}} - \varphi_{t}^{235 \overline{145}} - \varphi_{t}^{145 \overline{235}} \big).
        \end{split}
    \end{equation*}
    We have that $\widetilde{\Omega}_{0} = \Omega_{\ast}$ and from an easy calculation, we get that $d \widetilde{\Omega}_{t} = 0$, $\forall t \in I$. The transversality of $\widetilde{\Omega}_{t}$, $\forall t \in I$ can be proved in the same way as in Example \ref{etabeta_{5}}. Moreover, $\widetilde{\Omega}_{t} = e^{i_{\varphi(t)}|i_{\overline{\varphi(t)}}} (\Omega_{\ast})$.

    Thus, $\{\widetilde{\Omega}_{t}\}$ is a family of $3$-K\"ahler structures.      
\end{example}

\section{Cohomological aspects of p-K\"ahler manifolds}\label{Section 6}
In this section we study cohomological properties of $p$-K\"ahler, $p$-symplectic and $p$-pluriclosed structures.
\begin{theorem}\label{exactness of p-Kahler form}
    Let $(M,J,\Omega)$ be a compact $p$-K\"ahler manifold of complex dimension $n$, for \\
    $1 < p \leq n-1$. Let $k \doteq n-p$. If 
    \begin{equation*}
        [\Omega]_{\#} = 0, \quad \text{for $\# = dR, A, BC, \partial, \overline{\partial}$},  
    \end{equation*}
    then, there are no non-vanishing, simple, $(k,0)$-forms $\xi$ such that $\del \xi = \delbar \xi = 0$.
\end{theorem}
\begin{proof}
    Let us prove it for the Aeppli cohomology. Suppose that $[\Omega]_{A} = 0$, i.e., $\Omega = \del \lambda + \delbar \sigma$, for $\lambda \in \mathcal{A}^{p-1,p}(M), \sigma \in \mathcal{A}^{p,p-1}(M)$. If $\exists \, \xi \in \mathcal{A}^{k,0}(M)$ which is simple and such that $\del \xi = \delbar \xi = 0$, then 
    \begin{equation*}
        \begin{split}
            0 & < \int_{M} \sigma_{n-p} \, \Omega \wedge \xi \wedge \overline{\xi} = \int_{M} \sigma_{n-p} \big( \del \lambda + \delbar \sigma \big) \wedge \xi \wedge \overline{\xi} = \int_{M} \sigma_{n-p} \, \del \lambda \wedge \xi \wedge \overline{\xi} + \int_{M} \sigma_{n-p} \, \delbar \sigma \wedge \xi \wedge \overline{\xi} \\
            & = \int_{M} \sigma_{n-p} \, \del \big( \lambda \wedge \xi \wedge \overline{\xi} \big) + \int_{M} \sigma_{n-p} \, \delbar \big( \sigma \wedge \xi \wedge \overline{\xi} \big) \pm \int_{M} \sigma_{n-p} \, \lambda \wedge \big( \del \xi \big)\wedge \overline{\xi} \\
            & \quad \, \pm \int_{M} \sigma_{n-p} \, \lambda \wedge \xi \wedge \big(\del \overline{\xi} \big) \pm \int_{M} \sigma_{n-p} \, \sigma \wedge \big( \delbar \xi \big) \wedge \overline{\xi} \pm \int_{M} \sigma_{n-p} \, \sigma \wedge \xi \wedge \big( \delbar \overline{\xi} \big) \\
            & = \int_{M} \sigma_{n-p} \, d \big( \lambda \wedge \xi \wedge \overline{\xi} \big) + \int_{M} \sigma_{n-p} \, d \big( \sigma \wedge \xi \wedge \overline{\xi} \big) = 0.
        \end{split}
    \end{equation*}
    The proof for the other classes of cohomology follows in a similar way. 
\end{proof}
We can also prove a similar result for $p$-symplectic structures. 
\begin{theorem}\label{exactness of p-sympl}
    Let $(M,J,\Psi)$ be a compact $p$-symplectic manifold of complex dimension $n$, for $1 \leq p \leq n-1$. Let $k \doteq n-p$. If 
    \begin{equation*}
        [\Psi]_{dR} =0 ,  
    \end{equation*}
    then, there are no non-vanishing, simple, $(k,0)$-forms $\xi$ such that $\del \xi = \delbar \xi = 0$.
\end{theorem}
\begin{proof}
    Suppose that $[\Psi]_{dR} = 0$, i.e., $\Psi = d \eta$, for $\eta \in \mathcal{A}^{2p-1}(M)$. If $\exists \, \xi \in \mathcal{A}^{k,0}(M)$ which is simple and such that $\del \xi = \delbar \xi = 0$, then 
    \begin{equation*}
        0  <  \int_{M} \sigma_{n-p} \, \Psi^{p,p} \wedge \xi \wedge \overline{\xi} = \int_{M} \sigma_{n-p} \, \Psi \wedge \xi \wedge \overline{\xi} = \int_{M} \sigma_{n-p} \big( d \eta) \wedge \xi \wedge \overline{\xi},
    \end{equation*}
    where $\Psi^{p,p}$ is the $(p,p)$-component of $\Omega$. 

    Thus 
    \begin{equation*}
        0 < \int_{M} \sigma_{n-p} \big( d \eta  \big) \wedge \xi \wedge \overline{\xi} = \int_{M} \sigma_{n-p} \, d \big( \eta  \wedge \xi \wedge \overline{\xi} \big) \pm \int_{M} \sigma_{n-p} \, \eta  \wedge \big( d \xi \big) \wedge \overline{\xi} \pm \int_{M} \sigma_{n-p} \, \eta  \wedge \xi \wedge \big( d\overline{\xi} \big) = 0.
    \end{equation*}
\end{proof}

A similar result also holds for $p$-pluriclosed structures. 

\begin{theorem}\label{Necessary condition for the vanishing of p-pluriclosed forms}
    Let $(M,J,\Omega)$ be a compact $p$-pluriclosed manifold of complex dimension $n$, for $1 \leq p \leq n-1$. Let $k \doteq n-p$. If 
    \begin{equation*}
        [\Omega]_{A} =0 ,  
    \end{equation*}
    then, there are no non-vanishing, simple, $(k,0)$-forms $\xi$ such that $\del \xi = \delbar \xi = 0$.
\end{theorem}
\begin{proof}
    The proof is analogue to the proof of Theorem \ref{exactness of p-Kahler form}
\end{proof}

\begin{remark}
    We mention that, for $p=1,n-1$, Theorem \ref{Necessary condition for the vanishing of p-pluriclosed forms} coincides with \cite[Theorem 1.1]{PT2020}.
\end{remark}

As last application, we study the cohomology classes of $p$-K\"ahler, $p$-symplectic and $p$-pluriclosed structures on nilmanifolds equipped with invariant complex structures.

\medskip

Let $(M, J)$ be a nilmanifold of complex dimension $n$ equipped with an invariant complex structure $J$. According to \cite[Theorem 1.3]{S2001}, there exists $\{\varphi^{1}, \dots \varphi^{n}\}$ basis of left-invariant $(1,0)$-forms such that 
\begin{equation*}
    d \varphi^{j} \in \text{I} \, (\varphi^{1}, \dots, \varphi^{j-1}) \quad j = 1, \dots, n.
\end{equation*}
Therefore, $d \varphi^{1}=0$ and 
\begin{equation*}
    d (\varphi^{1} \wedge \dots \wedge \varphi^{j} ) = 0, \quad \text{for $j=2, \dots n$}.
\end{equation*}
Thus, we have the following theorem.
\begin{theorem}
    Let $(M,J)$ be a nilmanifold of complex dimension $n$. Let $1 \leq p \leq n-1$. Then
    \begin{enumerate}
        \item if $(M,J)$ admits a $p$-K\"ahler structure $\Omega$, then 
        \begin{equation*}
            0 \neq [\Omega]_{\#} \in H_{\#}(M), \quad \text{for $\#=dR,A,BC,\del, \delbar$} \,;
        \end{equation*}
        \item if $(M,J)$ admits a $p$-symplectic structure $\Psi$, then 
        \begin{equation*}
            0 \neq [\Psi]_{dR} \in H^{2p}_{dR}(M);
        \end{equation*}
        \item if $(M,J)$ admits a $p$-pluriclosed structure $\Omega$, then 
        \begin{equation*}
            0 \neq [\Omega]_{A} \in H^{p,p}_{A}(M).
        \end{equation*}
    \end{enumerate}
\end{theorem}

\end{document}